\documentclass{article}

\usepackage{xcolor}
\usepackage{parskip}
\usepackage{hyperref}
\usepackage{algorithm}
\usepackage[noend]{algorithmic}
\usepackage{amsmath,amssymb, amsthm}
\usepackage[a4paper, total={6in, 8in}]{geometry}
\usepackage[
    backend=biber,
    style=authoryear,
    citestyle=authoryear,
    natbib=true,
    maxnames=2,
    minnames=1,
    maxbibnames=99,
    dashed=false,
    uniquename=false,
    uniquelist=true,
    labeldateparts=true,
    doi=false,
    isbn=false,
    url=false,
    bibencoding=utf8
]{biblatex}

\AtBeginBibliography{\emergencystretch=3em}

\newcommand{\wrt}{wrt.\@}

\newcommand{\lhs}{left hand side}
\newcommand{\T}{^{\top}}

\newcommand{\iid}{i.i.d.\@}
\newcommand{\polylog}[1]{\operatorname{polylog}(#1)}

\newcommand{\infnorm}[1]{\|#1\|_{\infty}}

\newcommand{\onenorm}[1]{\|#1\|_{1}}
\newcommand{\twonorm}[1]{\|#1\|_{2}}

\newcommand{\abs}[1]{|#1|}

\newcommand{\Twonorm}[1]{\left\|#1\right\|_2}

\newcommand{\Abs}[1]{\left|#1\right|}

\newcommand{\E}{\mathbb{E}}
\newcommand{\var}[2][]{\operatorname{Var}_{#1}(#2)}

\newcommand{\indicator}[1]{\mathbf{1}{\{#1\}}}

\renewcommand{\P}{\mathbb{P}}

\newcommand{\cardinality}[1]{|#1|}
\newcommand{\R}{\mathbb{R}}
\newcommand{\N}{\mathbb{N}}

\DeclareMathOperator{\eigmax}{\lambda_{max}}

\newcommand{\proj}[2]{\Pi_{#1}(#2)}

\newcommand{\parent}[1]{(#1)}
\newcommand{\bracket}[1]{[#1]}
\newcommand{\curly}[1]{\{{#1}\}}

\newcommand{\Parent}[1]{\left(#1\right)}
\newcommand{\Bracket}[1]{\left[#1\right]}
\newcommand{\Curly}[1]{\left\{{#1}\right\}}

\theoremstyle{plain}

\newtheorem{definition}{Definition}[section]
\newtheorem{lemma}[definition]{Lemma}
\newtheorem{theorem}[definition]{Theorem}
\newtheorem{corollary}[definition]{Corollary}

\newtheorem{assumption}[definition]{Assumption}

\theoremstyle{definition} 

\newtheorem{example}[definition]{Example}

\newtheorem*{remark}{Remark}

\title{Simple Private Noisy SGD is Minimax-Optimal \\ under Population Strong Convexity}

\author{
  Valentin Roth\thanks{Department of Statistics, University of Oxford}
  \and
  Marco Avella Medina\thanks{Department of Statistics, Columbia University}
}

\date{October 2026}

\begin{document}

\maketitle

\begin{abstract}

    Stochastic gradient descent (SGD) is at the center of modern machine learning and is known to be statistically minimax-optimal in many regression settings, e.g., for convex Generalized Linear Models (GLMs). Yet, it is unknown whether its simplest differentially private Noisy SGD counterpart fulfills similar favorable guarantees even on strongly convex losses. We resolve this question affirmatively. We upper bound the parameter estimation error for a single pass of Noisy SGD with a fixed stepsize schedule and constant batch sizes on smooth and strongly convex Lipschitz losses over a compact convex domain under Renyi DP. Specialized to strongly convex GLMs, our error bounds match minimax lower bounds we derive based on the Score Attack framework. In particular, this provides the first minimax-optimal rate for Noisy SGD in private logistic regression, avoiding poly-logarithmic terms present in existing results. Our analysis is versatile as it relies on standard convex analysis and establishes the same rates in expectation and with high probability. 

\end{abstract}

\section{Introduction}
\label{sec:introduction}

Many supervised machine learning tasks can be formulated as stochastic convex optimization (SCO) problems. In regression, e.g., a population risk $\mathcal L_P(\theta) := \E\Bracket{\mathcal L_n(\theta)} := \E\Bracket{\frac 1 n \sum_{i=1}^n \ell(\theta, (Y_i, X_i))}$ is to be minimized, where $X_i$ are \iid{} observations from a prespecified distribution with true parameter $\vartheta \in \Theta \subseteq \R^d$ that shall be recovered. In times of big data, SCO is mostly solved using stochastic gradient-based optimization procedures where $\mathcal L_P$ is approximated using the loss $\ell$ evaluated on a subsample of the observations instead of the whole empirical risk $\mathcal L_n$. For classical regression models, stochastic gradient descent (SGD) can be shown to achieve statistically optimal performance when recovering $\vartheta$ through an SCO with properly specified loss components $\ell$ \citep{polyak:juditsky1992, bach:moulines2013, toulis:rennie:airoldi2014}. When data privacy guarantees like approximate, Gaussian or Renyi differential privacy (DP) \citep{dwork:mcsherry:nissim:smith2006, dong:roth:su2022, mironov:2017} need to be enforced, the picture is less clear even in simple convex and strongly convex settings. This work focuses on the simplest possible case, regression problems that can be solved through SCO with strongly convex and smooth Lipschitz losses with a single pass of SGD. 

Objective perturbation solves such problems by privatizing the empirical risk through noise and solving this new optimization problem using standard convex optimization methods. Similarly, output perturbation algorithms solve the non-private problem first, e.g., using SGD and afterwards privatize the estimate by adding noise. Both approaches are known to be optimal for recovering the empirical/population risk minimizer for strongly convex and smooth problems (See \cite{chaudhuri:monteleoni:sarwate2011, kifer:smith:thakurta2012, talwar:thakurta2014} for objective and \cite{zhang:zheng:mou:wang2017, lowy:razaviyayn2021} for output perturbation).

Noisy SGD-type algorithms follow another approach and privatize every single gradient within the iterations of SGD via additive Gaussian noise. This approach was introduced by \cite{song:chaudhuri:sarwate2013} and later popularized under the name of DP-SGD by \cite{abadi:chu:goodfellow:mcmahan2016}. Initially, the privacy of such algorithms was shown through composition theorems or privacy amplification by subsampling. Privacy guarantees based on amplification by subsampling analyze Noisy SGD's privacy when mini-batches are sampled uniformly at random, with or without replacement or using Poisson sampling \citep{abadi:chu:goodfellow:mcmahan2016, mironov:talwar:zhang2019, wang:balle:kasiviswanathan2019, dong:roth:su2022, bu:dong:long:weijie2020}. In a concurrent line of work, the privacy amplification by iteration guarantees in \cite{feldman:mironov:talwar:thakurta2018, feldman:koren:talwar2020} tighten privacy bounds by analyzing the Renyi DP leakage when only releasing the last iterate after one pass of Noisy SGD over the data, with fixed disjoint batches. Most recently, works by \cite{chourasia:ye:shokri2021, altschuler:talwar2022, ye:shokri2022} introduced convergent Renyi DP bounds, for full-batch and mini-batch Noisy SGD, where the privacy loss on (strongly) convex and smooth losses eventually stops growing after a fixed number of iterations. Likewise, \cite{bok:su:altschuler2024} proved this phenomenon under Gaussian DP. Such convergent Renyi DP bounds extend to non-convex and also non-smooth losses \citep{assodeh:diaz2023, altschuler:bok:talwar2024, kong:ribero2024, chien:li2025, liang:zhang:he:wu:xing2026}. Unfortunately, the generic convergent privacy bounds seem to be too weak for optimal utility guarantees under strong convexity. For instance, the privacy bounds requiring strong convexity like \cite[Theorem 5.1]{altschuler:bok:talwar2024} are stated for constant stepsizes, tie the privacy inducing noise to a lower bound on the batch size or introduce excess $\polylog{n}$ factors. Others like \cite[Theorem 4.2]{bok:su:altschuler2024} and \cite[Theorem 4.2]{ye:shokri2022} do not provide closed form expressions required for deriving precise guarantees. They also all impose strong convexity on the loss components $\ell$. However, in $\R^d$ with $d > 1$, the components $\ell$ used for parameter estimation for important regression models -- like generalized linear models (GLMs)\footnote{See, e.g., \cite{nelder:wedderburn1972} for an introduction.} -- are not strongly convex without regularization. In turn, the empirical risk $\mathcal L_n$ or population risk $\mathcal L_P$ may still be strongly convex.

As a consequence, the (almost) optimal utility guarantees for gradient-based algorithms in the literature rely on simpler privacy bounds and make up for this through algorithm design such as complicated batch size schedules etc. \cite{bassily:smith:thakurta2014} systematically study the optimization performance of several algorithms. Thereby, they provide an $(\varepsilon,\delta)$-DP optimization lower bound of order $\frac{d}{\varepsilon^2 n^2}$ \cite[Theorem 5.5]{bassily:smith:thakurta2014} and show that full-batch Noisy SGD achieves this rate up to $\polylog{n}$ in $n^2$ iterations \citep[Theorem 2.4]{bassily:smith:thakurta2014}. Ultimately, their guarantees transfer to statistical estimation in GLMs \citep[Theorem F.4]{bassily:smith:thakurta2014}. \cite{feldman:koren:talwar2020} shave off this remaining $\polylog{n}$ using only $n$ samples through a complicated algorithm that repeatedly solves auxiliary optimization problems \citep[Theorem 5.1]{feldman:koren:talwar2020}. Alternatively, they provide a version of Noisy SGD that performs a single pass over the observations with an elaborate batch sizing scheme \citep[Theorem 5.4]{feldman:koren:talwar2020}, whose privacy is governed by privacy amplification by iteration. Yet, this algorithm again suffers an excess $\polylog{n}$ factor. Likewise, the parameter estimation upper bounds for GLMs in \cite{cai:wang:zhang2021, cai:wang:zhang2025} suffer from such polylogarithmic factors. 

It is unclear whether a \emph{simple and practical} gradient-based algorithm can achieve statistical optimality under privacy. We resolve this question affirmatively, relying on \cite{feldman:koren:talwar2020}'s privacy amplification by iteration Renyi DP guarantee. Concretely, we show that a single pass of Noisy SGD with fixed batch sizes achieves minimax estimation rates for parameter estimation in strongly convex GLMs. We verify the optimality of our upper bound derived using simple convex analysis techniques through a matching lower bound obtained via the Score Attack framework \citep{cai:wang:zhang2025}. Our work thereby contributes to a growing body of literature that evaluates the statistical performance of (full batch) Noisy SGD on regression problems under privacy \citep{cai:wang:zhang2021, cai:wang:zhang2025, avella:bradshaw:loh2023, brown:dvijotham:evans:liu:smith:thakurta2024, marchis:loh2025, bombari:luo:seroussi:mondelli2025, lin:kolaczyk:smith:paquette2025}. To the best of our knowledge, the minimax-optimal bound for linear regression in \cite{bombari:luo:seroussi:mondelli2025} is the only other purely Noisy SGD-based utility guarantee for GLMs without excess polylogarithmic factors. Given that their proofs rely on deterministic equivalents in the high-dimensional proportional regime where $d/n\to c\in(0,\infty)$, we want to highlight the simplicity of our analysis. 

\section{Preliminaries}

\subsection{Renyi Differential Privacy}

The privacy notion we use is Renyi differential privacy (Renyi DP) introduced by \cite{mironov:2017}. 

\begin{definition}{(Definition 4, \cite{mironov:2017}).}
    \label{def:renyiDP}
    For measures $P$ and $Q$ with $P \ll Q$, define
    \begin{align*}
        D_\alpha (P||Q)
        := \frac{1}{\alpha-1} \log\Parent{\int \Parent{\frac{dP}{dQ}}^\alpha dQ}
        \quad \text{for $\alpha \in (1,\infty)$}. 
    \end{align*}
    When $\alpha \in \curly{1, \infty}$, define $D_\alpha (P||Q)$ by continuity. If $X \sim P$ and $Y \sim Q$, we write $D_\alpha(X||Y)$. Now, a randomized Algorithm $\mathcal A : \mathcal X^n \to \mathcal Y$ fulfills $(\alpha, \epsilon)$-Renyi DP, if for all $x,y \in \mathcal X^n$ with Hamming distance $d_h(x,y) \leq 1$ for fixed $\alpha \in [1,\infty]$ and $\epsilon \geq 0$, we have $D_\alpha(\mathcal A(x) || \mathcal A(y)) \leq \epsilon$.
    
\end{definition}

\subsection{Empirical Risk Minimization}

We approach SCO through empirical risk minimization (ERM), which approximates the expectation \wrt{} the data within the population risk using its empirical counterpart given $n$ \iid{} observations from the same distribution. We define the empirical and population risk: 

\begin{definition}
    Let $x \in \mathcal X^n$ be a fixed and $X \in \mathcal X^n$ a random dataset. For loss components $\ell : \Theta \times \mathcal X \to \R$, define the empirical risk $\mathcal L_n$ and population risk $\mathcal L_P$ as 
    \begin{align*}
        \mathcal L_n(\theta) 
        := \frac 1 n \sum_{i=1}^n \ell(\theta, x_i)
        \quad \text{and} \quad 
        \mathcal L_P(\theta) 
        := \E\Bracket{\mathcal L_n(\theta)}
        = \E\Bracket{\frac 1 n \sum_{i=1}^n \ell(\theta, X_i)}. 
    \end{align*}
    We call $\hat\theta \in \hat\Theta := \arg\min_{\theta \in \Theta} \mathcal L_n(\theta)$ and $\theta^\star \in \Theta^\star := \arg\min_{\theta \in \Theta} \mathcal L_P(\theta)$ an empirical and population risk minimizer, respectively. 
\end{definition}

The two main assumptions we impose on the losses are strong convexity and smoothness. 

\begin{definition}
    \label{def:strongconvexity}
    We call a differentiable convex function $f : \Theta \to \R$ 
    \begin{itemize}
        \setlength\itemsep{-0.3em}
        \item $\mu$-strongly convex if $f(\theta) \geq f(\vartheta) + \nabla_\theta f(\vartheta)\T(\theta-\vartheta) + \frac{\mu}{2} \twonorm{\theta-\vartheta}^2, \forall \theta, \vartheta \in \Theta$.
        \item $\nu$-smooth if $f(\theta) \leq f(\vartheta) + \nabla_\theta f(\vartheta)\T(\theta-\vartheta) + \frac{\nu}{2} \twonorm{\theta-\vartheta}^2, \forall \theta, \vartheta \in \Theta$.
    \end{itemize}
\end{definition}

\subsection{Noisy Stochastic Gradient Descent}

Stochastic gradient descent (SGD) approximates the gradient of the population risk by the gradient of the batch risk $\nabla_\theta \mathcal L_{B_k}(\theta)$ in Definition \ref{def:batchrisk} instead of the gradient of the empirical risk $\nabla_\theta \mathcal L_n(\theta)$. For large datasets, this may reduce each iteration's computational burden. 

\begin{definition}
    \label{def:batchrisk}
    For $x \in \mathcal X^n$ and $B_k \subseteq [n]$, define $\mathcal L_{B_k}(\theta) := \frac{1}{\cardinality{B_k}} \sum_{i \in B_k} \ell(\theta, x_i)$.
\end{definition}

To ensure data privacy, Noisy SGD in addition injects Gaussian noise at each iteration. The exact procedure can be found in Algorithm \ref{alg:noisySGD} below. 

\begin{algorithm}[H]
    \caption{NoisySGD($\ell, x, \{\eta_k\}_{k=1}^K, \{\sigma^2_k\}_{k=1}^K, \theta_0, \curly{B_k}_{k=1}^K$)}
    \label{alg:noisySGD}
    \begin{algorithmic}
        \REQUIRE $\ell : \Theta \times \mathcal X \to \R$, $x \in \mathcal X^n$, $\theta_0 \in \Theta$, $\{\eta_k\}_{k=1}^K, \{\sigma^2_k\}_{k=1}^K \in [0,\infty)^K$, $\{B_k\}_{k=1}^K \subseteq [n]$

        \FORALL{$k \in [K]$}
            \STATE $Z_k \gets Z \sim \mathcal N(0,\sigma_k^2 \cdot I_d)$
            \STATE $\theta_k \gets \proj{\Theta}{\theta_{k-1} - \eta_k \parent{\nabla_\theta \mathcal L_{B_k}(\theta_{k-1}) + Z_k}}$ 
        \ENDFOR
        \RETURN $\theta_K$
    \end{algorithmic}
\end{algorithm}

If the loss components $\ell$ have bounded gradient sensitivity as in Definition \ref{def:sensitivity} and are convex and smooth, Noisy SGD is differentially private in terms of Renyi DP. The formal guarantee is given in Theorem \ref{thm:privampbyiter} and holds for a single pass over the data when it is split into fixed and disjoint batches.

\begin{definition}
    \label{def:sensitivity}
    Let $x_i \sim x_i^\prime$ denote that $x_i, x_i^\prime \in \mathcal X$ and $d_H(x_i,x_i^\prime) \leq 1$. A loss component $\ell : \Theta \times \mathcal X \to \R$ has bounded gradient sensitivity $S_g$ if 
    \begin{align*}
        \sup_{\theta \in \Theta, x_i \sim x_i^\prime} \twonorm{\nabla_\theta \ell(\theta, x_i) - \nabla_\theta \ell(\theta, x_i^\prime)}
        \leq S_g
        < \infty.
    \end{align*}
\end{definition}

\begin{theorem}{(Theorem 3.1, \cite{feldman:koren:talwar2020}).}
    \label{thm:privampbyiter}
    Let $\Theta \subseteq \R^d$ be a convex set. Moreover, let $\ell$ be convex and $\nu_\ell$-smooth in $\theta \in \Theta$ and have bounded gradient sensitivity $S_g < \infty$.  Let $\{\eta_k\}_{k=1}^K, \{\sigma^2_k\}_{k=1}^K\in [0,\infty)^K$ be stepsize and noise schedules and $\curly{B_k}_{k=1}^K$ a schedule of disjoint batches s.t.\ $\bigsqcup_{k \in [K]} B_k = [n]$. Assume that $\eta_k \leq \frac{2}{\nu_\ell}$ for all $k \in [K]$ and let $\theta_K$ and $\theta^\prime_K$ be the outputs of Algorithm \ref{alg:noisySGD} initialized at $\theta_0 \in \Theta$ and run on $x, x^\prime \in \mathcal X^n$, respectively. Then, 
    \begin{align*}
        D_\alpha(\theta_K || \theta^\prime_K) 
        \leq \alpha S_g^2 \max_{k \in [K]} \Curly{\frac{\eta_k^2}{\cardinality{B_k}^2 \sum_{s=k}^K \eta_s^2 \sigma_s^2}}
        \quad \text{for all $\alpha \geq 1$}. 
    \end{align*}
\end{theorem}

\begin{remark}
    \label{rem:conversion}
    Assume $\frac{\eta_{k+1}}{\cardinality{B_{k+1}}} \leq \frac{\eta_k}{\cardinality{B_k}}$ for all $k \in [K]$. Similar to \cite[Equation (4.2)]{bombari:luo:seroussi:mondelli2025}, we observe that to achieve a $(\alpha,\epsilon)$-Renyi DP guarantee given a noise or stepsize schedule, respectively, the stepsizes and noise variances should be chosen as follows: 
    \begin{align*}
        \eta_k
        = \sqrt{\frac \epsilon {\alpha S_g^2}} \cdot \sqrt{\cardinality{B_k}^2 \sum_{s=k}^K \eta_s^2 \sigma_s^2}
        \quad \text{and} \quad 
        \sigma_k^2
        = \frac{\alpha S_g^2}{\epsilon} \begin{cases}
            \frac{1}{\cardinality{B_k}^2} - \frac{1}{\eta_k^2}\frac{\eta_{k+1}^2}{\cardinality{B_{k+1}}^2} & \text{if $k \in [K-1]$} \\
            \frac{1}{\cardinality{B_K}^2} & \text{if $k = K$}. 
        \end{cases}
    \end{align*}
\end{remark}

\section{Private Parameter Estimation}
\label{sec:parameterestimation}

We focus on the parameter estimation error of Noisy SGD and control $\twonorm{\theta_K-\theta^\star}$ both in expectation and with high probability, where $\theta_K$ is the output of Algorithm \ref{alg:noisySGD} and $\theta^\star$ is the population risk minimizer. This is in contrast to analyzing the optimization error in terms of $\twonorm{\theta_K-\hat\theta}$ or $\mathcal L_n(\theta_K)-\mathcal L_n(\hat\theta)$ as it directly incorporates the statistical error. This has multiple advantages. It avoids unnecessary triangle inequalities and union bounding and enables imposing a strong-convexity assumption directly onto the population risk $\mathcal L_P$ instead of the empirical risk $\mathcal L_n$. Establishing strong convexity for $\mathcal L_P$ is usually easier than for $\mathcal L_n$. 

\subsection{Setting up the Analysis} 

As is typical for the analysis of gradient-based optimization algorithms on strongly convex functions, we first establish that the parameter estimation error contracts for each iteration of Algorithm \ref{alg:noisySGD}. This is stated in Lemma \ref{lem:estimationerrorcontraction} in Appendix \ref{app:general}, which in addition to the standard contraction also contains the additive SGD noise term $\eta_k^2 \twonorm{\Delta_{1,k}}^2$ and privacy noise term $\eta_k^2 \twonorm{Z_k}^2$, as well as three cross terms. This is similar to \cite[Lemma 6]{chourasia:ye:shokri2021} but we contract around $\theta^\star$ as opposed to $\hat\theta$ and account for our mini-batch setting. Recursing Lemma \ref{lem:estimationerrorcontraction}, from the first iterate $\theta_1$ to the last one $\theta_K$, yields Lemma \ref{lem:estimationerror}, proven in Appendix \ref{app:general}. Its bound is random both in the data $X \in \mathcal X^n$ and the privacy noise $Z_1,...,Z_K$.
\begin{lemma}
    \label{lem:estimationerror}
    Let $X \in \mathcal X^n$ with \iid{} $X_i$ and $\Theta \subset \R^d$ be convex. Let $\mathcal L_P$ be $\mu$-strongly convex and $\nu$-smooth over $\Theta$ and $\eta_1,...,\eta_K \leq \frac{2}{\mu+\nu}$. Then, for the $K$-th iterate of Algorithm \ref{alg:noisySGD},
    \begin{align*}
        \twonorm{\theta_K-\theta^\star}^2 
        &\leq \underbrace{\twonorm{\theta_0-\theta^\star}^2 \cdot \Pi_{0,K}^{\mu,\nu}}_{=: T_{\operatorname{Bias}}} + \underbrace{\sum_{k=1}^K 2\eta_k^2 \twonorm{\Delta_{1,k}}^2 \cdot \Pi_{k,K}^{\mu,\nu}}_{=: T_{\operatorname{SGD}}} + \underbrace{\sum_{k=1}^K 2\eta_k^2 \twonorm{Z_k}^2 \cdot \Pi_{k,K}^{\mu,\nu}}_{=: T_{\operatorname{DP}}} + T_{\operatorname{Cross}}. 
    \end{align*}
    Here, the cross term is defined as 
    \begin{align*}
        T_{\operatorname{Cross}}
        :=& - \sum_{k=1}^K 2\eta_k (\theta_{k-1}-\theta^\star-\eta_k \Delta_{2,k})\T(\Delta_{1,k} + Z_k) \cdot \Pi_{k,K}^{\mu,\nu}
    \end{align*}
    and the quantities $\Pi_{k,K}^{\mu,\nu}$, $\Delta_{1,k}$ and $\Delta_{2,k}$ in the above are defined as follows: 
    \begin{align*}
        \begin{aligned}
            \Pi_{k,K}^{\mu,\nu} 
            := \prod_{j=k+1}^K \Parent{1-\frac{2\eta_j \mu\nu}{\mu+\nu}} \qquad
            \begin{aligned}
                \Delta_{1,k} &:= \nabla_\theta \mathcal L_{B_k}(\theta_{k-1}) - \nabla_\theta \mathcal L_P(\theta_{k-1}), \\
                \Delta_{2,k} &:= \nabla_\theta \mathcal L_P(\theta_{k-1}) - \nabla_\theta \mathcal L_P(\theta^\star).
            \end{aligned}
        \end{aligned}
    \end{align*}
\end{lemma}

While the cross term $T_{\operatorname{Cross}}$ in Lemma \ref{lem:estimationerror} vanishes in expectation, controlling it is challenging in probability (See Subsection \ref{sec:probability}). 

\subsection{Analysis in Expectation}
\label{sec:expectation}

Our in-expectation analysis builds on Lemma \ref{lem:estimationerror}. We first show that the cross term $\E[T_{\operatorname{Cross}}]$ indeed vanishes in expectation. To do so, we use iterated expectations conditioning on
\begin{align*}
    \mathcal F_{k-1}
    := \sigma\Parent{\theta_0, \curly{X_i : i \in B_1 \cup \dotsm \cup B_{k-1}}, Z_1,...,Z_{k-1}}.
\end{align*} 
Note that $\Delta_{2,k}$ is $\mathcal F_{k-1}$-measurable and thus $\E[T_{\operatorname{Cross}}] = 0$ if $\E\Bracket{\Delta_{1,k}|\mathcal F_{k-1}} = 0$ and $\E[Z_k|\mathcal F_{k-1}] = 0$. The former is true as $\nabla_\theta \mathcal L_{B_k}(\theta_{k-1})$ is unbiased for $\nabla_\theta \mathcal L_P(\theta_{k-1})$ when $\theta_{k-1}$ is fixed. The latter holds because $Z_k \sim \mathcal N(0, \sigma_k^2 I_d)$ is drawn independently. Combining $\E[T_{\operatorname{Cross}}] = 0$ with the fact that $\E[\twonorm{Z_k}^2|\mathcal F_{k-1}] = d\sigma^2$, we have essentially derived Lemma \ref{lem:esterrorinexp}. Its formal proof in Appendix \ref{app:expectation} shaves of the factor 2 from the terms  $T_{\operatorname{SGD}}$ and $T_{\operatorname{DP}}$ in Lemma \ref{lem:estimationerror}. This comes at the cost of controlling another cross term. 

\begin{lemma}
    \label{lem:esterrorinexp}
    Let $X \in \mathcal X^n \subseteq \R^{n \times d}$ with \iid{} $X_i$ and let $\Theta \subseteq \R^d$ be convex. Let $\mathcal L_P$ be $\mu$-strongly convex and $\nu$-smooth,  $\eta_1,...,\eta_K \leq \frac{2}{\mu+\nu}$ and $\curly{B_k}_{k=1}^K$ be disjoint. Then, 
    \begin{align*}
        \E\Bracket{\twonorm{\theta_K - \theta^\star}^2}
        &\leq \E\Bracket{\twonorm{\theta_0 - \theta^\star}^2} \cdot \Pi_{0,K}^{\mu,\nu} + \sum_{k=1}^K \eta_k^2 \Parent{\E\Bracket{\twonorm{\Delta_{1,k}}^2} + d\sigma_k^2} \cdot \Pi_{k,K}^{\mu,\nu}.
    \end{align*}
\end{lemma}

Given the conditional unbiasedness of $\nabla_\theta \mathcal L_{B_k}(\theta_{k-1})$, the terms $\E\bracket{\twonorm{\Delta_{1,k}}^2|\mathcal F_{k-1}}$ within $\E\Bracket{T_{\operatorname{SGD}}}$ are actually the variances of the stochastic gradient $\nabla_\theta \mathcal L_{B_k}(\theta_{k-1})$ when $\theta_{k-1}$ is fixed. A uniformly bounded gradient $\sup_{\theta \in \Theta, x_i \in \mathcal X} \twonorm{\nabla_\theta \ell(\theta,x_i)} \leq \frac{L}{2} < \infty$ implies an upper bound of order $\frac{L^2}{\cardinality{B_k}}$ on these variances (See Corollary \ref{cor:iidbatchvariance} in Appendix \ref{app:expectation}). Moreover, it implies a bounded gradient sensitivity through a triangle inequality. That makes the Renyi DP guarantee in Theorem \ref{thm:privampbyiter} applicable. Yet, the theorem's smoothness assumption is slightly stronger than the one in Lemma \ref{lem:esterrorinexp}. Theorem \ref{thm:privampbyiter} imposes $\nu_\ell$-smoothness on the loss components $\ell$ instead of $\nu$-smoothness on $\mathcal L_P$. Smoothness of $\mathcal L_P$ is then implied, because 
\begin{align*}
    \eigmax\Parent{\nabla_\theta^2 \mathcal L_P(\theta)}
    \leq \E\Bracket{\max_{\twonorm{u}=1} u\T \nabla_\theta^2 \ell(\theta, X_i) u}
    \leq \max_{x_i \in \mathcal X} \eigmax\Parent{\nabla_\theta^2 \ell(\theta,x_i)}
    \leq \nu_\ell. 
\end{align*}
\begin{remark}
    An immediate consequence is that $1 \leq \frac{\nu_\ell}{\mu}$ as $\mu I_d \preceq \nabla_\theta^2 \mathcal L_P(\theta) \preceq \nu_\ell I_d$. 
\end{remark}

To combine Theorem \ref{thm:privampbyiter} and Lemma \ref{lem:esterrorinexp}, it remains to specify a stepsize schedule $\{\eta_k\} \in [0,\infty)^K$. This schedule must both balance $\E[T_{\operatorname{Bias}}]$ and $\E[T_{\operatorname{SGD}}]$ as in prior non-private analyses of SGD and ensure that the corresponding noise schedule $\{\sigma^2_k\}_{k=1}^K \in [0,\infty)^K$ implied by Remark \ref{rem:conversion} makes the privacy term $\E[T_{\operatorname{DP}}]$ vanish at the correct rate. On top, Theorem \ref{thm:privampbyiter} and Lemma \ref{lem:esterrorinexp} require $\eta_k \leq \frac{2}{\mu+\nu_\ell} \leq \frac{2}{\nu_\ell}$. Stepsizes like $\eta_k = \frac{1}{\mu k}$ as in \cite{rakhlin:shamir:sridharan2011, harvey:liaw:plan:randhawa2019} or $\eta_k = \frac{2}{\mu(k+1)}$ as in \cite{harvey:liaw:randhawa2019} used for the analysis of non-private SGD on strongly convex and smooth losses fail this condition for general $\mu \leq \nu_\ell$. Additionally, they introduce excess polylogarithmic factors given the privacy term. We use the more general family $\eta_k = \frac{\mu+\nu_\ell}{2\mu\nu_\ell}\frac{\beta}{k+\tau}$, where $\beta, \tau \in \N$ remain to be specified. Choosing $\beta \geq 3$ removes polylogarithmic overhead and specifying the right $\tau$ ensures that $\mu \leq \nu_\ell$ can be arbitrary. We are ready to state our first main result as Theorem \ref{thm:esterrorinexp}. It provides an in-expectation error bound for the last iterate of noisy SGD for arbitrary partitions of $[n]$ of (nearly) uniform size $b$. 

\begin{theorem}
    \label{thm:esterrorinexp}
    Let $X \in \mathcal X^n \subseteq \R^{n \times d}$ with \iid{} $X_i \in \mathcal X$ and let $\Theta \subseteq \R^d$ be convex. Further, let $\mathcal L_P$ be $\mu$-strongly convex and let $\ell$ be convex and $\nu_\ell$-smooth in $\theta \in \Theta$, with $\sup_{\theta \in \Theta, x_i \in \mathcal X} \twonorm{\nabla_\theta \ell(\theta, x_i)} \leq \frac{L}{2} < \infty$. Let $K = \lfloor \frac{n}{b} \rfloor$ and $\cardinality{B_k} = b$ for all $k \in [K-1]$ and $\cardinality{B_K} \geq b$. Let $\eta_k = \frac{\mu+\nu_\ell}{2\mu\nu_\ell}\frac{\beta}{k+\tau} \leq \frac{2}{\mu+\nu_\ell}$ for all $k \in [K]$ where $\beta, \tau \in \N$ and $3 \leq \beta-1 \leq \tau < K$. Then, Algorithm \ref{alg:noisySGD} is $(\alpha, \epsilon)$-Renyi DP and 
    \begin{align*}
        \E\Bracket{\twonorm{\theta_K - \theta^\star}^2}
        &\leq \E\Bracket{\twonorm{\theta_0-\theta^\star}^2} \cdot \Parent{\frac{\tau}{\lfloor\frac{n}{b}\rfloor}}^\beta + \beta^2 \Parent{\frac{L^2}{\mu^2\lfloor\frac{n}{b}\rfloor b} + \frac{3\alpha dL^2}{\mu^2\epsilon b^2 \lfloor\frac{n}{b}\rfloor^2}}.
    \end{align*}
\end{theorem}

\begin{remark}
    \label{rem:schedule}
    We can always construct schedules $\{\eta_k\} \in [0,\infty)^K$ fulfilling the assumptions of Theorem \ref{thm:esterrorinexp}. Since $\frac{1}{\nu_\ell} \leq \frac{2}{\mu+\nu_\ell} \leq \frac{2}{\nu_\ell}$, the choices $\beta = 4$ and $\tau = \lceil \frac{\nu_\ell}{\mu} + 2\rceil$ ensure that 
    \begin{align*}
        \max_{k \in [n]} \eta_k 
        = \eta_1
        \leq \frac{2}{\mu+\nu_\ell}
        \; \Leftrightarrow \; 
        \frac{\mu+\nu_\ell}{2\mu\nu_\ell} \frac{\beta}{1+\tau} 
        \leq \frac{2}{\mu+\nu_\ell}
        \; \Leftarrow \; 
        \frac{\nu_\ell}{\mu} +2 \leq \tau. 
    \end{align*}
    Choosing $\beta = 4$ above is arbitrary and simply made for convenience. 
\end{remark}

\begin{remark} 
    \label{rem:rate}
    Note that $\Parent{\frac{\tau}{K}}^\beta \leq \frac{\tau}{K} \asymp \frac{\nu_\ell}{\mu K}$ as $\tau < K$. Further, if $n\mod b = 0$ the $\lfloor \cdot \rfloor$ simplify, e.g., $Kb = \lfloor \frac{n}{b} \rfloor b = n$ for $b=1$. Combining this with $\twonorm{\theta_0-\theta^\star}^2 \leq \frac{L^2}{\mu^2}$ which follows from Corollary \ref{cor:ubiniterror} and assuming that $\frac{\nu_\ell}{\mu} \asymp 1$, we see that Theorem \ref{thm:esterrorinexp} simplifies to
    \begin{align*}
        \E\Bracket{\twonorm{\theta_n - \theta^\star}^2}
        &\lesssim \frac{L^2}{\mu^2n} + \frac{\alpha dL^2}{\mu^2\epsilon n^2}.
    \end{align*}
    This matches the private parametric minimax rate $\frac{d}{n} + \frac{d^2}{\epsilon n^2}$ when $\frac{L^2}{\mu^2} \asymp d$, as is the case in bounded mean estimation. Beyond that, Section \ref{sec:minimax} derives minimax-optimality for strongly convex GLMs with Lipschitz negative log-likelihood losses.
\end{remark}

The theorem further suggests that choosing large constant batch sizes $b \asymp 1$ is not necessarily advantageous as it does not affect the constants in the second and third term -- the SGD and privacy noise terms -- but increases the constant in the bias term. 

\subsection{Analysis in Probability}
\label{sec:probability}

The analysis in probability is also based on Lemma \ref{lem:estimationerror} and resembles the one in expectation. The SGD noise term $T_{\operatorname{SGD}}$ is bounded the same way using uniform boundedness of the gradient. To avoid excess polylogarithmic factors due to union bounding, the privacy noise term is concentrated as a whole. For this, we write $Z_k = \sigma_k X_k$ for $X_k \sim \mathcal N(0,I_d)$, yielding
\begin{align*}
    T_{\operatorname{DP}}
    := \sum_{k=1}^K \eta_k^2 \twonorm{Z_k}^2 \cdot \Pi_{k,K}^{\mu,\nu_\ell}
    = \sum_{k=1}^K \sum_{j=1}^d \underbrace{\sigma_k^2 \eta_k^2 \Pi_{k,K}^{\mu,\nu_\ell}}_{=: w_{kj}} (X_k)_j^2.
\end{align*}
Using the Gaussian Hanson-Wright inequality for sums of weighted Gaussians (See Theorem \ref{thm:hansonwright} in Appendix \ref{app:probability}), we then obtain that with probability at least $1-\delta$,
\begin{align*}
    T_{\operatorname{DP}}
    \leq \Parent{1+\sqrt{\log(1/\delta)}}^2 d \sum_{k=1}^K \sigma_k^2 \eta_k^2 \Pi_{k,K}^{\mu,\nu_\ell}. 
\end{align*}
Unlike in expectation, the term $T_{\operatorname{Cross}}$ does not easily vanish here. While trivial bounds based on Cauchy Schwarz and Young's inequality do not recover the correct convergence rate, we obtain it using a generalized version of Freedman's inequality in Theorem \ref{thm:freedman}. In Remark \ref{rem:freedman} we elaborate that Theorem \ref{thm:freedman} is actually a direct corollary of \cite[Lemma C.3]{harvey:liaw:plan:randhawa2019}.

\begin{theorem}{(Generalized Freedman; Lemma C.3, \cite{harvey:liaw:plan:randhawa2019}).}
    \label{thm:freedman}
    Let $\Curly{D_t}_{t=1}^T$ be a martingale difference sequence adapted to the filtration $\Curly{\mathcal F_t}_{t=1}^T$. For all $t \in [T]$, let $v_{t-1} \geq 0$ be $\mathcal F_{t-1}$-measurable and suppose that the $D_t$ are conditionally $v_{t-1}$-sub-Gaussian:
    \begin{align*}
        \E\bracket{\exp(\lambda D_t)|\mathcal F_{t-1}} \leq \exp\Parent{\frac{\lambda^2 v_{t-1}}{2}} 
        \quad \text{for all $\lambda > 0$}. 
    \end{align*}
    If $\sum_{t=1}^T v_{t-1} \leq \sum_{t=1}^T \alpha_t D_t + \beta(\delta)$ with probability at least $1-\delta$ and $\alpha_1,...,\alpha_T,\beta(\delta) \geq 0$, 
    \begin{align*}
        \P\Parent{\sum_{t=1}^T D_t \lesssim \infnorm{\alpha} \log\Parent{1/\delta} + \sqrt{\beta(\delta) \log\Parent{1/\delta}}}
        \geq 1-2\delta. 
    \end{align*}
\end{theorem}

The theorem introduced in \cite{harvey:liaw:plan:randhawa2019} for the analysis of non-private SGD requires the sum of the conditional sub-Gaussian variance proxies to be bounded by an affine transformation of the martingale differences. The martingale differences $D_k$ in the cross term
\begin{align*}
    T_{\operatorname{Cross}}
    = - \sum_{k=1}^K 2\eta_k (\theta_{k-1}-\theta^\star-\eta_k \Delta_{2,k})\T(\Delta_{1,k} + Z_k) \cdot \Pi_{k,K}^{\mu,\nu}
    =: \sum_{k=1}^K D_k
\end{align*}
have zero mean, when conditioning on $\mathcal F_{k-1}$ as discussed previously in Subsection \ref{sec:expectation}. They are also adapted to $\curly{\mathcal F_k}_{k=1}^K$. Thus, they are martingale differences. Because they consist of bounded and Gaussian terms, they are also conditionally sub-Gaussian with variance proxy
\begin{align*}
    v_{k-1}
    \leq 4\eta_k^2 \Parent{\Pi_{k,K}^{\mu,\nu_\ell}}^2 (G^2+\sigma_k^2) \cdot \Parent{1-\frac{2\eta_k\mu\nu}{\mu+\nu}}\twonorm{\theta_{k-1}-\theta^\star}^2.
\end{align*}
Since $\twonorm{\theta_{k-1}-\theta^\star}^2$ appears, applying Lemma \ref{lem:estimationerror} introduces the martingale differences $D_k$ into the bound on $\sum_{t=1}^T v_{t-1}$. To apply the theorem it only remains to separate the coefficients $\infnorm{\alpha}$ and $\sqrt{\beta(\delta)}$. They turn out to be bounded up to constants by $T_{\operatorname{Bias}}, T_{\operatorname{SGD}}$ and $T_{\operatorname{DP}}$. Thus, $T_{\operatorname{Cross}}$ is of the same order and we obtain Lemma \ref{lem:esterrorinprobgen}.

\begin{lemma}
    \label{lem:esterrorinprobgen}
    Let $\Theta \subseteq \R^d$ be convex and $X \in \mathcal X^n \subseteq \R^{n \times d}$ with \iid{} $X_i$. Let $\mathcal L_P$ be $\mu$-strongly convex and $\nu$-smooth and $\sup_{k \in [K]} \twonorm{\Delta_{1,k}} \leq G < \infty$. Then, with probability at least $1-2\delta$ where $\delta \in (0, \frac{1}{2})$, Algorithm \ref{alg:noisySGD} with $\eta_1,...,\eta_K \leq \frac{2}{\mu+\nu}$ and disjoint $\curly{B_k}_{k=1}^K$ fulfills
    \begin{align*}
        \twonorm{\theta_K-\theta^\star}^2
        &\lesssim \Parent{\twonorm{\theta_0-\theta^\star}^2 \cdot \Pi_{0,K}^{\mu,\nu} + G^2 \sum_{k=1}^K \eta_k^2 \Pi_{k,K}^{\mu,\nu} + d\log(1/\delta) \sum_{k=1}^K \sigma_k^2 \eta_k^2 \Pi_{k,K}^{\mu,\nu}} \log(1/\delta). 
    \end{align*}
\end{lemma}

The proof of Lemma \ref{lem:esterrorinprobgen} in Appendix \ref{app:probability} requires one additional assumption over Lemma \ref{lem:esterrorinexp}, its in-expectation analogue. We assume $\twonorm{\Delta_{1,k}} \leq G < \infty$, which immediately follows from a uniformly bounded gradient with $G = L$. Instantiating Lemma \ref{lem:esterrorinprobgen} with $b=1$ and the stepsizes discussed in Remark \ref{rem:schedule} yields a high-probability analogue of Theorem \ref{thm:esterrorinexp}. 

\begin{theorem}
    \label{thm:esterrorinprob}
    Let $X \in \mathcal X^n \subseteq \R^{n \times d}$ with \iid{} $X_i$ and let $\Theta \subseteq \R^d$ be convex. Let $\mathcal L_P$ be $\mu$-strongly convex and let $\ell$ be convex and $\nu_\ell$-smooth in $\theta \in \Theta$, with $\sup_{\theta \in \Theta, x_i \in \mathcal X} \twonorm{\nabla_\theta \ell(\theta, x_i)} \leq \frac{L}{2} < \infty$. Suppose that $b = 1$, $\eta_k = \frac{\mu+\nu_\ell}{2\mu\nu_\ell}\frac{4}{k+\lceil\nu_\ell/\mu+2\rceil}$ for all $k \in [n]$ and $\lceil \frac{\nu_\ell}{\mu}+2\rceil \leq n$. Then, Algorithm \ref{alg:noisySGD} is $(\alpha, \epsilon)$-Renyi DP and with probability at least $1-2\delta$
    \begin{align*}
        \twonorm{\theta_n-\theta^\star}^2
        \lesssim \Parent{\twonorm{\theta_0-\theta^\star}^2 \cdot \frac{\nu_\ell}{\mu n} + \frac{L^2}{\mu^2 n} + \frac{\alpha dL^2}{\mu^2 \epsilon n^2}\log(1/\delta)}\log(1/\delta).
    \end{align*}
\end{theorem}

Theorem \ref{thm:esterrorinprob} is proven in Appendix \ref{app:probability}. Its high-probability upper bound is structurally analogous to the in-expectation bound in Theorem \ref{thm:esterrorinexp} besides the $\log(1/\delta)$ factors stemming from the high-probability control on the privacy noise and martingale differences. This is particularly remarkable as no additional polylogarithmic factors in $n$ are introduced. Hence, the minimax rates discussed next are also achievable in probability. 

\section{Minimaxity}
\label{sec:minimax}

To argue that the rates above are statistically optimal, we derive a private minimax lower bound and matching upper bound for a class of GLMs, including the logistic model. For that, we define a non-private and private class of estimators $\mathcal M$ and $\mathcal M_{\alpha,\epsilon}$ in Definition \ref{def:classes}. 

\begin{definition}
    \label{def:classes}
    We define the following classes of estimators $M : \mathcal Z^n \to \Theta \subseteq \R^d$: 
    \begin{itemize}
        \setlength\itemsep{-0.3em}
        \item $\mathcal M$ is the class of all measurable estimators $M$, 
        \item $\mathcal M_{\alpha,\epsilon}$ is the class of all $(\alpha,\epsilon)$-Renyi DP estimators $M \in \mathcal M$. 
    \end{itemize}
\end{definition}

\subsection{Score Attack}

We prove our minimax lower bound via the Score Attack framework \citep{cai:wang:zhang2025}. A Score Attack as in Definition \ref{def:score} is a generalization of fingerprinting codes and tracing attacks \citep{gabor2008,bun:ullman:vadhan2014,dwork:smith:steinke:ullman:vadhan2015, cai:wang:zhang2021,kamath:mouzakis:singhal2022}. It is constructed from the score statistic and governed by the Fisher information. We provide a new score attack result in Lemma \ref{lem:ubexpscoreattack} below which adapts \cite[Proposition 2.1]{cai:wang:zhang2025} given under $(\varepsilon,\delta)$-DP to Renyi DP. We note that \cite{bombari:luo:seroussi:mondelli2025} used a similar approach to derive minimax lower bounds for zero-concentrated DP estimators on linear models. 

\begin{definition}
    \label{def:score}
    Let $\Theta \subseteq \R^d$ and $\mathcal P_\Theta := \curly{P_\theta : \theta \in \Theta}$ be a parametric family of distributions with densities $p_\theta : \mathcal Z \to (0,\infty)$. Define the score statistic and Fisher information as 
    \begin{align*}
        S_\theta(z) := \nabla_\theta \log\Parent{p_\theta(z)}
        \quad \text{and} \quad 
        I(\theta) 
        := \var{S_\theta(Z)}
        \quad \text{where $Z \sim P_\theta$ and $\theta \in \Theta$}. 
    \end{align*}
    For $M \in \mathcal M_{\alpha,\epsilon}$ and $z \in \mathcal Z$, define the Score Attack as $\mathcal A_\theta\Parent{z, M(Z)} := \Parent{M(Z)-\theta}\T S_\theta(z)$. 
\end{definition}

The Score Attack has a hypothesis testing interpretation and is said to be \emph{sound} if it is small if the candidate $z$ is not part of the training dataset $Z$ used by the estimator $M(Z)$. It is \emph{complete} if the attack is large when $z$ is contained in $Z$ \citep{cai:wang:zhang2025}.  Theorem \ref{thm:soundness} shows soundness as it states that the attack evaluated on $z$ not in $Z$ has mean zero and variance bounded by the product of the maximum eigenvalue of $I(\theta)$ and the mean squared error (MSE) of $M(Z)$. 

\begin{theorem}{(Theorem 2.1, \cite{cai:wang:zhang2025}).}
    \label{thm:soundness}
    Let $Z := (Z_1,...,Z_n)\T \in \mathcal Z^n$ where $Z_i \sim P_\theta$ with $\theta \in \Theta \subseteq \R^d$. For all $i \in [n]$, suppose $Z_{(i)} \in \mathcal Z^n$ is an independent copy of $Z$ where the $i$-th entry is replaced by $Z_i^\prime \sim P_\theta$. Further, assume that $I(\theta) < \infty$, $\E\Bracket{S_\theta(Z)} = 0$ and $\E\Bracket{\twonorm{M(Z) - \theta}^2} < \infty$ where $Z \sim P_\theta$. Then, for all $i \in [n]$,
    \begin{align*}
        \E\Bracket{\mathcal A_\theta\Parent{Z_i, M(Z_{(i)})}} = 0
        \quad \text{and} \quad 
        \E\Bracket{\mathcal A_\theta\Parent{Z_i, M(Z_{(i)})}^2}
        \leq \eigmax\Parent{I(\theta)} \E\Bracket{\twonorm{M(Z)-\theta}^2}. 
    \end{align*}
\end{theorem}

\begin{remark}
    The second part of Theorem \ref{thm:soundness} is slightly more general than in \cite[Theorem 2.1]{cai:wang:zhang2025}. It is proven as an intermediate result in the proof of their Theorem 2.1.  
\end{remark}

For Renyi DP estimators soundness can be transferred to the Score Attack even when evaluated on $z = Z_i$ which is part of the training dataset $Z$. This is because by the definition of Renyi DP the distributions of $M(Z)$ and $M(Z_{(i)})$ cannot differ a lot. 

\begin{lemma}
    \label{lem:ubexpscoreattack}
    Let $M \in \mathcal M_{\alpha,\epsilon}$ with $\alpha \geq 2$ and $Z := (Z_1,...,Z_n)\T \in \mathcal Z^n$ where $Z_i \sim P_\theta$ \iid{} with $\theta \in \Theta \subseteq \R^d$. For all $i \in [n]$, let $Z_{(i)} \in \mathcal Z^n$ be an independent copy of $Z$ where the $i$-th entry is replaced by $Z_i^\prime \sim P_\theta$. If $I(\theta) < \infty$, $\E\Bracket{S_\theta(Z_i)} = 0$ and $\E\Bracket{\twonorm{M(Z) - \theta}^2} < \infty$,
    \begin{align*}
        \E\Bracket{\mathcal A_\theta\Parent{Z_i,M(Z)}}
        &\leq \sqrt{e^\epsilon-1} \cdot \sqrt{\eigmax\Parent{I(\theta)} \cdot \E\Bracket{\twonorm{M(Z)-\theta}^2}}. 
    \end{align*}
\end{lemma}

The proof of Lemma \ref{lem:ubexpscoreattack} in Appendix \ref{app:minimax} exploits that the likelihood ratio of Renyi DP algorithms trained on two neighboring datasets is controlled by the Renyi divergence. 

\subsection{Generalized Linear Model}

We specialize the parametric family to Generalized Linear Models (GLMs) as defined in Definition \ref{def:GLM}. The family of interest to us is then $\mathcal P_{\Theta_0} := \curly{\operatorname{GLM}_{\Theta_0}(\theta) : \theta \in \Theta_0}$. 

\begin{definition}
    \label{def:GLM}
    Fix a domain $\Theta_0 := \curly{\theta \in \R^d : \twonorm{\theta} \leq R_\Theta < \infty}$, partition function $\psi : \R \to \R$, base density $h : \mathcal Y \subseteq \R \times (0,\infty) \to [0,\infty)$ and covariate distribution $P_{\mathcal X}$ on $\R^d$. $(Y_1,X_1) \in \mathcal Z := \mathcal Y \times \mathcal X \subseteq \R^{d+1}$ is generated by $\operatorname{GLM}_{\Theta}^\psi(\vartheta)$ if $Y_1 \in \mathcal Y$ has conditional density
    \begin{align*}
        p_{\vartheta} (y|x)
        := \exp\Parent{\frac{y x\T\vartheta - \psi(x\T\vartheta)}{c(\varsigma)}} h(y,\varsigma). 
    \end{align*}
\end{definition}

\subsubsection{Statistical Lower Bound}

To prove the minimax lower bound, we instantiate the upper bound on the Score Attack in Lemma \ref{lem:ubexpscoreattack} for GLMs, sum the upper bound over the observations $Z_1,...,Z_n$ and take expectation \wrt{} a specific prior $\pi$ over $\vartheta$. This yields the first inequality below:  
\begin{align*}
    n \sqrt{e^\epsilon-1} \sqrt{\E\Bracket{\E\Bracket{\twonorm{M((Y,X))-\vartheta}^2|\vartheta}}}
    \gtrsim \sum_{i=1}^n \E\Bracket{\E\Bracket{\mathcal A_\vartheta\Parent{(Y_i,X_i), M((Y,X))}\Big|\vartheta}}
    \gtrsim d. 
\end{align*}
The second inequality is a `prior-averaged completeness bound` due to \cite{cai:wang:zhang2025} using the same prior (See Lemma \ref{lem:lbexpscoreattack} in Appendix \ref{app:minimax}). Solving for the expected MSE on the \lhs{} and noticing that the supremum over $\vartheta \in \Theta_0$ upper bounds the expectation over the prior, we obtain Lemma \ref{lem:lbpriv} whose full proof can be found in Appendix \ref{app:minimax}. 

\begin{lemma}
    \label{lem:lbpriv}
    Let $(Y, X) \in \mathcal Z^n \subseteq \R^{n \times (d+1)}$ contain \iid{} samples $(Y_i, X_i\T)\T \in \mathcal Z \subseteq \R^{d+1}$ from $\operatorname{GLM}_{\Theta_0}^\psi(\vartheta)$. If $\alpha \geq 2$, $\frac{d}{n} \leq \sqrt{\epsilon}$ and assuming $\infnorm{\psi^{\prime\prime}} \leq C_\psi^{\prime\prime}$ and $\E\Bracket{X_iX_i\T} \preceq \xi I_d \prec \infty$,
    \begin{align*}
        \inf_{M \in \mathcal M_{\alpha,\epsilon}}\sup_{\vartheta \in \Theta_0} \E\Bracket{\twonorm{M((Y,X))-\vartheta}^2|\vartheta} 
        \gtrsim \Parent{R_\Theta^2 \wedge \frac{c(\varsigma)}{C_\psi^{\prime\prime} \xi}} \cdot \frac{d^2}{(e^\epsilon-1) n^2}. 
    \end{align*}
\end{lemma}

For non-private estimators $M \in \mathcal M$, a minimax lower bound of order $\frac{d}{n}$ is well-established in the literature and can be obtained through Assouad's Lemma (See e.g., \cite[Section 8]{duchi:2023}). For completeness, we state it in Lemma \ref{lem:lbnonpriv} in Appendix \ref{app:minimax}. Combining this lower bound with Lemma \ref{lem:lbpriv} and using that for bounded $\epsilon \lesssim 1$ we have $e^\epsilon-1 \lesssim \epsilon$, we obtain Theorem \ref{thm:lbcombined}. Note that the theorem explicitly considers estimators on a convex and compact parameter space $\Theta_0 \subseteq \R^d$, matching the setting of our in-expectation upper bound in Theorem \ref{thm:esterrorinexp}.

\begin{theorem}
    \label{thm:lbcombined}
    Let $(Y, X) \in \mathcal Z \subseteq \R^{n \times (d+1)}$ contain \iid{} samples $(Y_i, X_i\T)\T \in \mathcal Z \subseteq \R^{d+1}$ from $\operatorname{GLM}_{\Theta_0}^\psi(\vartheta)$. Assuming $\alpha \geq 2$, $\frac{d}{n} \leq \sqrt{\epsilon} \lesssim 1$, $\infnorm{\psi^{\prime\prime}} \leq C_\psi^{\prime\prime}$ and $\E\Bracket{X_iX_i\T} \preceq \xi I_d \prec \infty$,
    \begin{align*}
        \inf_{M \in \mathcal M_{\alpha,\epsilon}}\sup_{\vartheta \in \Theta_0} \E\Bracket{\twonorm{M((Y,X))-\vartheta}^2|\vartheta} 
        \gtrsim \Parent{R_\Theta^2 \wedge \frac{c(\varsigma)}{C_\psi^{\prime\prime} \xi}} \cdot \Parent{\frac{d}{n} + \frac{d^2}{\epsilon n^2}}. 
    \end{align*}
\end{theorem}

\subsubsection{Statistical Upper Bound}

Finally, we also specialize Theorem \ref{thm:esterrorinexp} to GLMs that fulfill Assumption \ref{ass:covariates} and \ref{ass:partition}. We will see that when $\alpha \asymp 1$, the resulting upper bound matches the lower bound in Theorem \ref{thm:lbcombined}.  

\begin{assumption}
    \label{ass:covariates}
    For all $i \in [n]$, suppose that the covariates $X_i \sim P_{\mathcal X}$ on $\R^d$are a zero-mean $\upsilon^2$-sub-Gaussian vector. Furthermore, $\E\Bracket{X_iX_i\T} \succeq \zeta I_d \succ 0$, and $\twonorm{X_i} \leq \frac{R_{\mathcal X}}{2}\sqrt{d}$. 
\end{assumption}

\begin{assumption}
    \label{ass:partition}
    Assume that $\abs{y_1} \leq R_{\mathcal Y}$ for all $y_1 \in \mathcal Y$ and $\psi : \R \to \R$ is s.t.\ $\infnorm{\psi^\prime} \leq C_\psi^\prime$ and $\infnorm{\psi^{\prime\prime}} \leq C_\psi^{\prime\prime}$, $\psi^{\prime\prime}$ is symmetric, positive and $u \mapsto \psi^{\prime\prime}(\abs{u})$ is non-increasing. 
\end{assumption}

The loss components used to prove Lemma \ref{lem:esterrorinexpGLMs} in Appendix \ref{app:minimax} are the negative log-likelihood loss components $\ell(\theta, Z_i) := -\log\Parent{p_\theta(Y_i|X_i)}$ corresponding to the likelihood $\mathfrak L(\theta) =\prod_{i=1}^n p_\theta(Y_i,X_i) \propto \prod_{i=1}^n p_\theta(Y_i|X_i)$. Here, proportionality is up to constants in $\theta$. 

\begin{lemma}
    \label{lem:esterrorinexpGLMs}
    Let $(Y, X) \in \mathcal Z^n \subseteq \R^{n \times (d+1)}$ contain \iid{} samples $(Y_i, X_i\T)\T \in \mathcal Z \subseteq \R^{d+1}$ from $\operatorname{GLM}_{\Theta_0}^\psi(\vartheta)$. Let $P_{\mathcal X}$ and $\psi$ fulfill Assumption \ref{ass:covariates} and \ref{ass:partition}, respectively. Then, Algorithm \ref{alg:noisySGD} on $\ell(\theta, Z_i) := -\log\Parent{p_\theta(Y_i|X_i)}$ is $(\alpha,\epsilon)$-Renyi DP and with an appropriate $\{\eta_k\} \in [0,\infty)^n$, 
    \begin{align*}
        \E\Bracket{\twonorm{\theta_n - \vartheta}^2} 
        &\lesssim \frac{R_\Theta^2 R_{\mathcal X}^2}{\zeta \psi^{\prime\prime}(M_{\upsilon,\zeta,\Theta})} \Parent{\frac{d}{n}} + \frac{R_{\mathcal Y}^2 R_{\mathcal X}^2}{\zeta^2 \psi^{\prime\prime}(M_{\upsilon,\zeta,\Theta})^2} \Parent{\frac{d}{n} + \frac{\alpha d^2}{\epsilon n^2}}.
    \end{align*}
    This bound on the parameter estimation error requires that $n \gtrsim \frac{R_{\mathcal X}^2 d}{\zeta \psi^{\prime\prime}(M_{\upsilon,\zeta,\Theta})}$. 
\end{lemma}

\begin{example}
    Arguably the most important GLM with $\psi$ fulfilling Assumption \ref{ass:partition} is the logistic model. There, $c(\varsigma) = 1$, $\mathcal Y = \curly{0,1}$ $\psi(u) := \log(1+e^u)$, $\psi^\prime(u) = \sigma(u)$ and $\psi^{\prime\prime}(u) = \sigma^\prime(u)$ where $\sigma(u) = \frac{1}{1+e^{-u}}$ is the sigmoid function. Since $\infnorm{\sigma} \leq 1$, $\infnorm{\sigma^{\prime}} \leq \frac{1}{4}$ and $R_{\mathcal Y} = 1$, Assumption \ref{ass:partition} is fulfilled. For $P_{\mathcal X}$ s.t.\ Assumption \ref{ass:covariates} holds, by Corollary \ref{cor:sigmaprimelb} we have
    \begin{align*}
        \psi^{\prime\prime}(M_{\upsilon,\zeta,\Theta}) 
        = \sigma^\prime(M_{\upsilon,\zeta,\Theta}) 
        \geq \frac{1}{4} e^{-M_{\upsilon,\zeta,\Theta}}
        = \frac{1}{4} e^{-\tilde O(\upsilon R_\Theta)}. 
    \end{align*}
    For $n \gtrsim \frac{R_{\mathcal X}^2 d}{\zeta} e^{\tilde O(\upsilon R_\Theta)}$, the resulting upper bound on the parameter estimation error is 
    \begin{align*}
        \E\Bracket{\twonorm{\theta_n - \vartheta}^2} 
        &\lesssim e^{\tilde O(\upsilon R_\Theta)} \cdot R_{\mathcal X}^2 \cdot \Parent{\frac{R_\Theta^2}{\zeta} \Parent{\frac{d}{n}} + \frac{1}{\zeta^2} \Parent{\frac{d}{n} + \frac{\alpha d^2}{\epsilon n^2}}}.
    \end{align*}

    Compared to the non-private rate in \cite[Example 3.6]{puchkin:zhivotovskiy2023}, our bound only suffers exponentially in $\upsilon$ instead of the radius of $\mathcal X$, which here is of order $R_{\mathcal X} \sqrt{d}$. Yet, without privacy the exponential dependence on $R_\Theta$ can be reduced to $e^{O(\twonorm{\vartheta})}$ and removed fully if stronger assumptions are made \citep{chardon:lerasle:mourtada2024}. Whether this is possible under privacy remains open. 
\end{example}

Other GLMs fulfilling Assumption \ref{ass:partition} are binomial and Rademacher models. It remains to be investigated whether a comparably simple analysis extends to more general GLMs, e.g., through regularization and clipping to ensure strong convexity and Lipschitzness. 

\printbibliography

\appendix
\section{Proofs of Section \ref{sec:parameterestimation}}

The following identities are direct consequences of the dominated convergence theorem.
\begin{corollary}
    \label{cor:regularity}
    Let $\ell : \Theta \times \mathcal X \to \R$ be twice differentiable and $X \in \mathcal X^n$ be random. Then, 
    \begin{itemize}
        \setlength\itemsep{-0.3em}
        \item $\nabla_\theta \mathcal L_P(\theta) 
        = \E\Bracket{\nabla_\theta \mathcal L_n(\theta)}$ if $\sup_{\theta \in \Theta, x_i \in \mathcal X} \twonorm{\nabla_\theta \ell(\theta, x_i)} \leq \frac{L}{2} < \infty$, 
        \item  $\nabla_\theta^2 \mathcal L_P(\theta)
        = \E\Bracket{\nabla_\theta^2 \mathcal L_n(\theta)}$ if in addition $\ell$ is $\nu_\ell$-smooth. 
    \end{itemize}
\end{corollary}

\subsection{Setting up the Analysis}
\label{app:general}

\begin{lemma}{(Theorem 2.1.12, \cite{nesterov2018}).}
    \label{lem:stronglyconvexcocoercivity}
    Let $f : \mathcal Y \subseteq \R^d \to \R$ be a twice differentiable $\mu$-strongly convex and $\nu$-smooth function. Then, for all $x,y \in \mathcal Y$ it holds that
    \begin{align*}
        \Parent{\nabla f(x) - \nabla f(y)}\T \Parent{x-y} \geq \frac{\mu \nu}{\mu + \nu} \twonorm{x-y}^2 + \frac{1}{\mu+\nu} \twonorm{\nabla f(x) - \nabla f(y)}^2. 
    \end{align*}
\end{lemma}

\begin{lemma}{(Theorem 7.2 + Proposition 7.8, \cite{wright:recht2022}).}
    \label{lem:fixedpointoptimality}
    Let $\mathcal Y \subset \R^d$ be closed and convex, $f : \mathcal Y \to \R$ be differentiable and $\hat y \in \arg\min_{y \in \mathcal Y} f(y)$. For $\eta_k \geq 0$,
    \begin{align*}
        \hat y = \proj{\mathcal Y}{\hat y - \eta_k \nabla_y f(\hat y)}. 
    \end{align*}
\end{lemma}

\begin{corollary}
    \label{cor:paramerrorgradient}
    For a $\mu$-strongly convex function $f : \mathcal Y \subseteq \R^d \to \R$ and $\hat y \in \arg\min_{y\in \mathcal Y} f(y)$,
    \begin{align*}
        \twonorm{\hat y - y}
        \leq \frac{2}{\mu} \twonorm{\nabla_y f(y)}
        \quad \text{for all $y \in \mathcal Y$}. 
    \end{align*}
\end{corollary}

\begin{proof}
    By $\mu$-strong convexity of $f$ and since $f(\hat y) \leq f(y)$ for all $y \in \mathcal Y$, it holds that
    \begin{align*}
        f(y) 
        \geq f(\hat y) 
        \geq f(y) + \nabla_y f(y)\T (\hat y-y) + \frac{\mu}{2} \twonorm{\hat y - y}^2. 
    \end{align*}
    Subtracting $f(y)$ from both sides, rearranging and applying Cauchy Schwarz yields 
    \begin{align*}
        \twonorm{\hat y-y}^2
        \leq \frac{2}{\mu} \nabla_y f(y)\T (y-\hat y)
        \leq \frac{2}{\mu} \twonorm{\nabla_y f(y)} \twonorm{y-\hat y}. 
    \end{align*}
    If $\twonorm{y-\hat y} \neq 0$, dividing the above by that quantity, we have shown that
    \begin{align*}
        \twonorm{\hat y-y} \leq \frac{2}{\mu} \twonorm{\nabla_y f(y)}. 
    \end{align*}
    In the case when $\twonorm{y-\hat y} = 0$, the inequality above holds trivially.  
\end{proof}

\begin{lemma}
    \label{lem:estimationerrorcontraction}
    Let $\Theta \subset \R^d$ be a convex set and $X \in \mathcal X^n$ with \iid{} $X_i$. Let $\mathcal L_P$ be $\mu$-strongly convex and $\nu$-smooth over $\Theta$ and $\eta_1,...,\eta_k \leq \frac{2}{\mu+\nu}$. Then, for the $k$-th iterate of Algorithm \ref{alg:noisySGD}, 
    \begin{align*}
        \twonorm{\theta_k - \theta^\star}^2 
        &\leq \Parent{1-\frac{2\eta_k\mu\nu}{\mu+\nu}} \twonorm{\theta_{k-1}-\theta^\star}^2 + \eta_k^2 \Parent{\twonorm{\Delta_{1,k}}^2 + \twonorm{Z_k}^2} \\ & \quad - 2\eta_k \parent{\theta_{k-1}-\theta^\star-\eta_k\Delta_{2,k}}\T \Parent{\Delta_{1,k} + Z_k}  \\ & \quad + 2\eta_k^2 \Delta_{1,k} \T Z_k.
    \end{align*}
    The quantities $\Delta_{1,k}$ and $\Delta_{2,k}$ in the statement above are defined as follows: 
    \begin{align*}
        \underbrace{\nabla_\theta \mathcal L_{B_k}(\theta_{k-1}) - \nabla_\theta \mathcal L_P(\theta^\star)}_{=: \Delta_k}
        = \underbrace{\nabla_\theta \mathcal L_{B_k}(\theta_{k-1}) - \nabla_\theta \mathcal L_P(\theta_{k-1})}_{=:\Delta_{1,k}} + \underbrace{\nabla_\theta \mathcal L_P(\theta_{k-1}) - \nabla_\theta \mathcal L_P(\theta^\star)}_{=:\Delta_{2,k}}. 
    \end{align*}
\end{lemma}

\begin{proof}
    $\theta^\star = \proj{\Theta}{\theta^\star - \eta_k \nabla_\theta \mathcal L_P(\theta^\star)}$ by constrained first-order optimality of $\theta^\star$ for $\mathcal L_P$ in Lemma \ref{lem:fixedpointoptimality}. Combining this with the Noisy SGD update rule and non-expansiveness, 
    \begin{align*}
        \twonorm{\theta_k - \theta^\star}^2 
        &= \twonorm{\proj{\Theta}{\theta_{k-1} - \eta_k (\nabla_\theta \mathcal L_{B_k}(\theta_{k-1}) + Z_k)} - \theta^\star}^2 \\
        &= \twonorm{\proj{\Theta}{\theta_{k-1} - \eta_k (\nabla_\theta \mathcal L_{B_k}(\theta_{k-1}) + Z_k)} - \proj{\Theta}{\theta^\star - \eta_k \nabla_\theta \mathcal L_P(\theta^\star)}}^2 \\
        &\leq \twonorm{\theta_{k-1} - \eta_k (\nabla_\theta \mathcal L_{B_k}(\theta_{k-1}) + Z_k) - (\theta^\star - \eta_k \nabla_\theta \mathcal L_P(\theta^\star))}^2 \\
        &= \twonorm{\theta_{k-1} - \theta^\star - \eta_k (\underbrace{\nabla_\theta \mathcal L_{B_k}(\theta_{k-1}) - \nabla_\theta \mathcal L_P(\theta^\star)}_{= \Delta_k}) - \eta_k Z_k}^2. 
    \end{align*}
    By expanding the square we then obtain the following: 
    \begin{equation}
     \label{eq:update}  
     \begin{aligned}
     \twonorm{\theta_k - \theta^\star}^2 
        &\leq \twonorm{\theta_{k-1} - \theta^\star - \eta_k \Delta_k - \eta_k Z_k}^2 \\
        &= \twonorm{\theta_{k-1}-\theta^\star - \eta_k\Delta_k}^2 + \eta_k^2 \twonorm{Z_k}^2 - 2\eta_k \parent{\theta_{k-1}-\theta^\star-\eta_k\Delta_k}\T Z_k \\
        &= \twonorm{\theta_{k-1}-\theta^\star - \eta_k\Delta_{2,k}}^2 + \eta_k^2 \Parent{\twonorm{\Delta_{1,k}}^2 + \twonorm{Z_k}^2} \\ & \quad - 2\eta_k \parent{\theta_{k-1}-\theta^\star-\eta_k\Delta_{2,k}}\T \Delta_{1,k} \\ & \quad - 2\eta_k \parent{\theta_{k-1}-\theta^\star-\eta_k\Delta_{2,k}}\T Z_k \\ & \quad + 2\eta_k^2 \Delta_{1,k}\T Z_k.
        \end{aligned}
    \end{equation}
    We now show that the first term, i.e., the deterministic update, contracts: 
    \begin{align}
        \notag\twonorm{\theta_{k-1} &- \theta^\star - \eta_k \Delta_{2,k}}^2
        = \twonorm{\theta_{k-1}-\theta^\star}^2 + \eta_k^2 \twonorm{\Delta_{2,k}}^2 - 2\eta_k \parent{\theta_{k-1}-\theta^\star}\T \Delta_{2,k} \\
        \notag &\overset{(i)}\leq \twonorm{\theta_{k-1}-\theta^\star}^2 + \eta_k^2 \twonorm{\Delta_{2,k}}^2
        - 2\eta_k \Parent{\frac{\mu\nu}{\mu+\nu} \twonorm{\theta_{k-1}-\theta^\star}^2 + \frac{1}{\mu+\nu} \twonorm{\Delta_{2,k}}^2} \\
        \notag &= \Parent{1-\frac{\mu\nu}{\mu+\nu} 2\eta_k} \twonorm{\theta_{k-1}-\theta^\star}^2 + \eta_k \Parent{\eta_k - \frac{2}{\mu+\nu}} \twonorm{\Delta_{2,k}}^2 \\
        &\overset{(ii)}\leq \Parent{1-\frac{\mu\nu}{\mu+\nu} 2\eta_k} \twonorm{\theta_{k-1}-\theta^\star}^2. 
        \label{eq:deterministic_update}
    \end{align}
    Here, the inequality $(i)$ holds by Lemma \ref{lem:stronglyconvexcocoercivity} and $(ii)$ because $\eta_k \leq \frac{2}{\mu+\nu}$.
    Combining the statements in Displays \eqref{eq:update} and \eqref{eq:deterministic_update} yields the result. 
\end{proof}

\subsubsection{Proof of Lemma \ref{lem:estimationerror}}

\begin{proof}
    Recursing Lemma \ref{lem:estimationerrorcontraction}, as formalized in Lemma \ref{lem:unrollrec} yields 
    \begin{align*}
        \twonorm{\theta_K-\theta^\star}^2 
        &\leq \twonorm{\theta_0-\theta^\star}^2 \cdot \Pi_{0,K}^{\mu,\nu} + \sum_{k=1}^K \eta_k^2\Parent{\twonorm{\Delta_{1,k}}^2 + \twonorm{Z_k}^2} \cdot \Pi_{k,K}^{\mu,\nu} \\ & \quad\quad - \sum_{k=1}^K 2\eta_k (\theta_{k-1}-\theta^\star-\eta_k \Delta_{2,k})\T(\Delta_{1,k} + Z_k) \cdot \Pi_{k,K}^{\mu,\nu} \\ & \quad\quad + \sum_{k=1}^K 2\eta_k^2 \Delta_{1,k} \T Z_k \cdot \Pi_{k,K}^{\mu,\nu}. 
    \end{align*}
    Here, $\Delta_{1,k}, \Delta_{2,k}$ are as defined in Lemma \ref{lem:unrollrec}. At the cost of constants, the last sum above can be integrated into the first one. Concretely, by Cauchy Schwarz and Young's inequality, 
    \begin{align*}
        2\Delta_{1,k}\T Z_k 
        \leq 2 \twonorm{\Delta_{1,k}} \twonorm{Z_k}
        \leq \twonorm{\Delta_{1,k}}^2 + \twonorm{Z_k}^2. 
    \end{align*}
    Plugging this in recovers the statement.
\end{proof}

\subsection{Analysis in Expectation}
\label{app:expectation}

\begin{corollary}
    \label{cor:iidbatchriskisunbiased}
    Let $X \in \mathcal X^n$ with identically distributed $X_i$. Let $\sup_{\theta \in \Theta, x_i \in \mathcal X} \twonorm{\nabla_\theta \ell(\theta,x_i)} \leq \frac{L}{2} < \infty$. For the population and batch risk $\mathcal L_P$ and $\mathcal L_{B_k}$ with $B_k \subseteq [n]$, we have
    \begin{align*}
        \E\Bracket{\nabla_\theta \mathcal L_{B_k}(\theta)}
        = \E\Bracket{\nabla_\theta \ell(\theta,X_i)}
        = \nabla_\theta \mathcal L_P(\theta)
        \quad \text{for all $\theta \in \Theta \subseteq \R^d$}.
    \end{align*}
\end{corollary}

\begin{proof}
    By linearity of the expectation it holds that
    \begin{align*}
        \E\Bracket{\nabla_\theta \mathcal L_{B_k}(\theta)}
        &= \E\Bracket{\frac{1}{\cardinality{B_k}}\sum_{i \in B_k} \nabla_\theta \ell(\theta,X_i)}
        = \frac{1}{\cardinality{B_k}} \sum_{i \in B_k} \E\Bracket{\nabla_\theta \ell(\theta,X_i)} \\
        &\overset{(i)}= \frac{1}{n} \sum_{i=1}^n \E\Bracket{\nabla_\theta \ell(\theta,X_i)}
        = \E\Bracket{\nabla_\theta \mathcal L_n(\theta)} \\
        &\overset{(ii)}= \nabla_\theta \E\Bracket{\mathcal L_n(\theta)}
        = \nabla_\theta \mathcal L_P(\theta). 
    \end{align*}
    Here, $(i)$ uses the identical distribution and $(ii)$ holds by regularity shown in Corollary \ref{cor:regularity}. 
\end{proof}

\begin{corollary}
    \label{cor:iidbatchvariance}
    Let $X \in \mathcal X^n$ with \iid{} $X_i$ and $B_k \subseteq [n]$. Let the gradient be uniformly bounded, i.e., $\sup_{\theta \in \Theta, x_i \in \mathcal X} \twonorm{\nabla_\theta \ell(\theta,x_i)} \leq \frac{L}{2} < \infty$. Then, the variance of $\mathcal L_{B_k}$ is bounded as
    \begin{align*}
        \E\Bracket{\twonorm{\nabla_\theta \mathcal L_{B_k}(\theta) - \nabla_\theta \mathcal L_P(\theta)}^2}
        \leq \frac{L^2}{\cardinality{B_k}}
        \quad \text{for all $\theta \in \Theta \subseteq \R^d$}.
    \end{align*}
\end{corollary}

\begin{proof}
    For convenience, define $\delta_i := \nabla_\theta \ell(\theta,X_i) - \nabla_\theta \mathcal L_P(\theta)$ for all $\theta \in \Theta$. With this notation and using that the $X_i$ are \iid{} and that $\E\Bracket{\delta_i} = 0$ by Corollary \ref{cor:iidbatchriskisunbiased}, we get
    \begin{align*}
        \E\Bracket{\twonorm{\nabla_\theta \mathcal L_{B_k}(\theta) - \nabla_\theta \mathcal L_P(\theta)}^2}
        &= \E\Bracket{\Twonorm{\frac{1}{\cardinality{B_k}} \sum_{i \in B_k} \nabla_\theta \ell(\theta, X_i) - \nabla_\theta \mathcal L_P(\theta)}^2} \\
        &= \frac{1}{\cardinality{B_k}^2} \E\Bracket{\Twonorm{\sum_{i \in B_k} \delta_i}^2}\\
        &= \frac{1}{\cardinality{B_k}^2} \sum_{i \in B_k} \E\Bracket{\twonorm{\delta_i}^2} \\
        &\leq \frac{1}{\cardinality{B_k}^2} \sum_{i \in B_k} L^2 \\
        &= \frac{L^2}{\cardinality{B_k}}.
    \end{align*}
    The inequality above follows from
    \begin{align*}
        \E\Bracket{\twonorm{\delta_i}^2}
        = \E\Bracket{\twonorm{\nabla_\theta \ell(\theta,X_i)-\nabla_\theta \mathcal L_P(\theta)}^2}
        \leq 2 \E\Bracket{\Twonorm{\nabla_\theta \ell(\theta, X_i)}^2 + \twonorm{\nabla_\theta \mathcal L_P(\theta)}^2}
        \leq L^2. 
    \end{align*}
\end{proof}

\begin{corollary}{(See Lemma 5, \cite{rakhlin:shamir:sridharan2011}).}
    \label{cor:ubiniterror}
    For a $\mu$-strongly convex population risk $\mathcal L_P$ with uniformly bounded gradient, i.e., $\sup_{x_i \in \mathcal X} \twonorm{\nabla_\theta \ell(\theta_0, x_i)} \leq \frac{L}{2} < \infty$ we have 
    \begin{align*}
        \twonorm{\theta_0 - \theta^\star}
        \leq \frac{L}{\mu}. 
    \end{align*}
\end{corollary}

\begin{proof}
    The assumption of a uniformly bounded gradient, i.e., $\sup_{x_i \in \mathcal X} \twonorm{\nabla_\theta \ell(\theta_0, x_i)} \leq \frac{L}{2}$, Corollary \ref{cor:paramerrorgradient} and convex Jensen's applied to the norm $x \mapsto \twonorm{x}$ yield
    \begin{align*}
        \twonorm{\theta_0 - \theta^\star}
        &\leq \frac{2}{\mu} \twonorm{\nabla_\theta \mathcal L_P(\theta_0)}
        = \frac{2}{\mu} \twonorm{\E\bracket{\nabla_\theta \mathcal L_n(\theta_0)}} \\
        &\leq \frac{2}{\mu} \E\Bracket{\twonorm{\nabla_\theta \mathcal L_n(\theta_0)}}
        = \frac{2}{\mu n} \E\Bracket{\Twonorm{\sum_{i=1}^n \nabla_\theta \ell(\theta_0,x_i)}} \\
        &\leq \frac{2}{\mu n} \sum_{i=1}^n \E\bracket{\twonorm{\nabla_\theta \ell(\theta_0,x_i)}}
        \leq \frac{L}{\mu}. 
    \end{align*}
\end{proof}

\subsubsection{Proof of Lemma \ref{lem:esterrorinexp}}

\begin{proof}
    Using the intermediate result in the proof of Lemma \ref{lem:estimationerror}, taking expectation and applying the law of iterated expectations, we get
    \begin{align*}
        \E\Bracket{\twonorm{\theta_K-\theta^\star}^2}
        &\leq \E\Bracket{\twonorm{\theta_0-\theta^\star}^2} \cdot \Pi_{0,K}^{\mu,\nu} + \sum_{k=1}^K \eta_k^2 \E\Bracket{\E\Bracket{\twonorm{\Delta_{1,k}}^2 + \twonorm{Z_k}^2|\mathcal F_{k-1}}} \cdot \Pi_{k,K}^{\mu,\nu} \\ & \quad\quad - \sum_{k=1}^K 2\eta_k \E\Bracket{(\theta_{k-1}-\theta^\star-\eta_k \Delta_{2,k})\T \E\Bracket{\Delta_{1,k} + Z_k|\mathcal F_{k-1}}} \cdot \Pi_{k,K}^{\mu,\nu} \\ & \quad\quad + \sum_{k=1}^K 2\eta_k^2 \E\Bracket{\E\Bracket{\Delta_{1,k} \T Z_k|\mathcal F_{k-1}}} \cdot \Pi_{k,K}^{\mu,\nu}
    \end{align*}
    Above, we use that $\Delta_{2,k}$ is measurable \wrt{} the filtration $\curly{\mathcal F_k}_{k=1}^K$ where 
    \begin{align*}
        \mathcal F_k
        = \sigma\Parent{\theta_0, \curly{X_i : i \in B_1 \cup \dotsm \cup B_k}, Z_1,...,Z_k}.
    \end{align*}
    We now show that the two cross terms vanish, starting with the first one. 

    It holds that $\E\Bracket{Z_k|\mathcal F_{k-1}} = 0$ because $Z_k \sim \mathcal N(0, \sigma_k^2 I_d)$ and it is independent of $\theta_{k-1}$. Moreover, as the batches $\curly{B_k}_{k\in [K]}$ are disjoint, $\theta_{k-1}$ and the $\curly{X_i}_{i \in B_k}$ are independent, we may apply Corollary \ref{cor:iidbatchriskisunbiased} conditionally. This unbiasedness of the stochastic gradient yields 
    \begin{align*}
        \E\Bracket{\Delta_{1,k} + Z_k|\mathcal F_{k-1}}
        &= \E\Bracket{\Delta_{1,k}|\mathcal F_{k-1}} \\
        &= \E\Bracket{\nabla_\theta \mathcal L_{B_k}(\theta_{k-1}) - \nabla_\theta \mathcal L_P(\theta_{k-1})|\mathcal F_{k-1}} \\
        &= \E\Bracket{\nabla_\theta \mathcal L_{B_k}(\theta_{k-1})|\mathcal F_{k-1}} - \nabla_\theta \mathcal L_P(\theta_{k-1})
        = 0. 
    \end{align*}
    As $Z_k$, $\theta_{k-1}$ and $\curly{X_i}_{i \in B_k}$ are mutually independent and since $Z_k \sim \mathcal N(0,\sigma_k^2)$, 
    \begin{align*}
        \E\Bracket{\Delta_{1,k}\T Z_k|\mathcal F_{k-1}}
        = \E\Bracket{\Delta_{1,k}|\mathcal F_{k-1}}\T \E\Bracket{Z_k|\mathcal F_{k-1}}
        = 0. 
    \end{align*}
    As both are true for all $k \in [K]$, the cross terms are zero and we recover the statement. 
\end{proof}

\subsubsection{Proof of Theorem \ref{thm:esterrorinexp}}

\begin{proof}
    First, note that $\eta_k \leq \frac{2}{\mu+\nu_\ell} \leq \frac{2}{\nu_\ell}$ as $\mu \geq 0$.

    To show that Algorithm \ref{alg:noisySGD} is $(\alpha, \epsilon)$-Renyi DP, we bound the gradient sensitivity
    \begin{align*}
        \sup_{\theta \in \Theta, x_i, x_i^\prime \in \mathcal X} \Twonorm{\nabla_\theta \ell(\theta, x_i) - \nabla_\theta \ell(\theta, x_i^\prime)}
        \leq 2\sup_{\theta \in \Theta, x_i \in \mathcal X} \Twonorm{\nabla_\theta \ell(\theta, x_i)}
        \leq L.
    \end{align*}
    Combining the bounded gradient sensitivity over the convex $\Theta$ with the  $\nu_\ell$-smoothness of $\ell$ and given $\eta_k \leq \frac{2}{\nu_\ell}$, Theorem \ref{thm:privampbyiter} and Remark \ref{rem:conversion} ensure that Algorithm 1 is $(\alpha,\epsilon)$-Renyi DP. 

    Similarly, for utility, the fact that $\eta_k \leq \frac{2}{\mu+\nu_\ell}$ let's us instantiate Lemma \ref{lem:esterrorinexp}, which yields 
    \begin{align*}
        \E\Bracket{\twonorm{\theta_K - \theta^\star}^2}
        &\leq \underbrace{\E\Bracket{\twonorm{\theta_0 - \theta^\star}^2} \cdot \Pi_{0,K}^{\mu,\nu_\ell}}_{=: T_{\operatorname{Bias}}} + \underbrace{\sum_{k=1}^K \eta_k^2 \E\Bracket{\twonorm{\Delta_{1,k}}^2} \cdot \Pi_{k,K}^{\mu,\nu_\ell}}_{=: T_{\operatorname{SGD}}} + \underbrace{\sum_{k=1}^K \eta_k^2 d\sigma_k^2 \cdot \Pi_{k,K}^{\mu,\nu_\ell}}_{=: T_{\operatorname{DP}}}.
    \end{align*}
    Combined with Corollary \ref{cor:prodBias} this justifies that 
    \begin{align*}
        T_{\operatorname{Bias}}
        = \E\Bracket{\twonorm{\theta_0 - \theta^\star}^2} \cdot \Pi_{0,K}^{\mu,\nu_\ell}
        \leq \E\Bracket{\twonorm{\theta_0 - \theta^\star}^2} \cdot \Parent{\frac{\tau}{K}}^{\beta}.
    \end{align*}
    By Corollary \ref{cor:iidbatchvariance} and using that $\cardinality{B_k} \geq b$ for all $k \in [K]$, for all $k \in [K]$
    \begin{align*}
        \E\Bracket{\twonorm{\Delta_{1,k}}^2}
        &= \E\Bracket{\E\Bracket{\twonorm{\nabla_\theta \mathcal L_{B_k}(\theta_{k-1})-\nabla_\theta \mathcal L_P(\theta_{k-1})}^2|\mathcal F_{k-1}}}
        \leq \frac{L^2}{\cardinality{B_k}}
        \leq \frac{L^2}{b}. 
    \end{align*}
    Here, the filtration $\curly{\mathcal F_k}_{k=1}^K$ is defined through 
    \begin{align*}
        \mathcal F_k
        = \sigma\Parent{\theta_0, \curly{X_i : i \in B_1 \cup \dotsm \cup B_k}, Z_1,...,Z_k}.
    \end{align*}
    Combining the last inequality with Corollary \ref{cor:sumSGD} and $\frac{\mu+\nu_\ell}{2\mu\nu_\ell} \leq \frac{1}{\mu}$, we have
    \begin{align*}
        T_{\operatorname{SGD}}
        = \sum_{k=1}^K \eta_k^2 \E\Bracket{\twonorm{\Delta_{1,k}}^2} \cdot \Pi_{k,K}^{\mu,\nu_\ell}
        \leq \frac{L^2}{b} \sum_{k=1}^K \eta_k^2 \Pi_{k,K}^{\mu,\nu_\ell}
        \leq \frac{L^2}{b} \cdot \frac{\beta^2}{K} \Parent{\frac{\mu+\nu_\ell}{2\mu\nu_\ell}}^2
        \leq \frac{L^2\beta^2}{\mu^2 Kb}. 
    \end{align*}    
    Lastly, since $\cardinality{B_k} = b$ for all $k \in [K-1]$ and $\cardinality{B_K} \geq b$, by Corollary \ref{cor:sumDP}, it holds that
    \begin{align*}
        T_{\operatorname{DP}}
        = d\sum_{k=1}^K \sigma_k^2 \eta_k^2 \Pi_{k,K}^{\mu,\nu_\ell}
        \leq \frac{3\alpha dL^2 \beta^2}{\epsilon b^2 K^2} \Parent{\frac{\mu+\nu_\ell}{2\mu\nu_\ell}}^2
        \leq \frac{3\alpha dL^2 \beta^2}{\mu^2 \epsilon b^2 K^2}. 
    \end{align*}
    By collecting terms and plugging in $K = \lfloor \frac{n}{b} \rfloor$, we have shown the statement: 
    \begin{align*}
        \E\Bracket{\twonorm{\theta_K - \theta^\star}^2}
        &\leq T_{\operatorname{Bias}} + T_{\operatorname{SGD}} + T_{\operatorname{DP}} \\
        &\leq \E\Bracket{\twonorm{\theta_0-\theta^\star}^2} \cdot \Parent{\frac{\tau}{\lfloor\frac{n}{b}\rfloor}}^\beta + \frac{L^2\beta^2}{\mu^2 \lfloor\frac{n}{b}\rfloor b} + \frac{3\alpha dL^2\beta^2}{\mu^2\epsilon b^2 \lfloor\frac{n}{b}\rfloor^2}.
    \end{align*}
\end{proof}

\subsection{Analysis in Probability}
\label{app:probability}
\begin{remark}{(Generalized Freedman's Inequality).}
    \label{rem:freedman}
    \cite[Lemma C.3]{harvey:liaw:plan:randhawa2019} states 
    \begin{align*}
        \P\Parent{\sum_{t=1}^T D_t \geq x}
        \leq \delta + \exp\Parent{-\frac{x^2}{4x \cdot \infnorm{\alpha} + 8\beta(\delta)}}. 
    \end{align*}
    We simplify the above using the following equivalence and choose the positive root of $x$: 
    \begin{align*}
        \exp\Parent{-\frac{x^2}{4x \cdot \infnorm{\alpha} + 8\beta(\delta)}}
        \leq \delta
        \quad \Leftrightarrow \quad 
        x^2 - 4 \infnorm{\alpha} \log\Parent{1/\delta} x - 8\beta(\delta)\log\Parent{1/\delta} \geq 0. 
    \end{align*}
    The version in Theorem \ref{thm:freedman}, is then obtained using $\sqrt{a+b} \leq \sqrt{a} + \sqrt{b}$ for all $a, b \geq 0$. 
\end{remark}

\begin{theorem}{(Hanson-Wright Inequality; Lemma 1, \cite{laurent:massart2000}).}
    \label{thm:hansonwright}
    For a Gaussian random vector $X \sim \mathcal N(0, I_d)$ and $w \in \R^d$ with $w_1,...,w_d \geq 0$, we have
    \begin{align*}
        \P\Parent{\sum_{i=1}^d w_i (X_i^2-1) \geq 2 \twonorm{w} \sqrt{x} + 2\infnorm{w} x}
        \leq \exp\parent{-x}.
    \end{align*}
\end{theorem}

\subsubsection{Proof of Lemma \ref{lem:esterrorinprobgen}}

\begin{proof}
    Writing $M_k := -(\theta_{k-1}-\theta^\star-\eta_k \Delta_{2,k})\T(\Delta_{1,k} + Z_k)$, Lemma \ref{lem:estimationerror} yields 
    \begin{align*}
        \twonorm{\theta_K-\theta^\star}^2 
        &\leq \overbrace{\twonorm{\theta_0-\theta^\star}^2 \cdot \Pi_{0,K}^{\mu,\nu}}^{=: T_{\operatorname{Bias}}} + \overbrace{\sum_{k=1}^K 2\eta_k^2 \twonorm{\Delta_{1,k}}^2 \cdot \Pi_{k,K}^{\mu,\nu}}^{=: T_{\operatorname{SGD}}} \\ &\quad + \underbrace{\sum_{k=1}^K 2\eta_k^2 \twonorm{Z_k}^2 \cdot \Pi_{k,K}^{\mu,\nu}}_{=: T_{\operatorname{DP}}} + \underbrace{\sum_{k=1}^K 2\eta_k M_k \cdot \Pi_{k,K}^{\mu,\nu}}_{=: T_{\operatorname{Cross}}}.
    \end{align*}
    For the SGD noise term we use that $\twonorm{\Delta_{1,k}} \leq G$ and hence
    \begin{align*}
        T_{\operatorname{SGD}}
        = \sum_{k=1}^K 2\eta_k^2 \twonorm{\Delta_{1,k}}^2 \cdot \Pi_{k,K}^{\mu,\nu}
        \leq 2G^2 \sum_{k=1}^K \eta_k^2 \Pi_{k,K}^{\mu,\nu}.
    \end{align*}
    We have $Z_k \sim \mathcal N(0, \sigma_k^2I_d)$ and $Z_k = \sigma_k X_k$ for $X_k \sim \mathcal N(0,I_d)$ for all $k \in [K]$. Thus, 
    \begin{align*}
        T_{\operatorname{DP}}
        &= \sum_{k=1}^K 2\eta_k^2 \twonorm{Z_k}^2 \cdot \Pi_{k,K}^{\mu,\nu}
        = 2\sum_{k=1}^K \sigma_k^2 \eta_k^2 \Pi_{k,K}^{\mu,\nu} \twonorm{X_k}^2 \\
        &= 2\sum_{k=1}^K \sigma_k^2 \eta_k^2 \Pi_{k,K}^{\mu,\nu} \sum_{j=1}^d (X_k)_j^2
        = 2\sum_{k=1}^K \sum_{j=1}^d \underbrace{\sigma_k^2 \eta_k^2 \Pi_{k,K}^{\mu,\nu}}_{=: w_{kj}} (X_k)_j^2.
    \end{align*}
    This shows that $T_{\operatorname{DP}}$ is just a weighted sum of $Kd$ independent squared standard Gaussians. We collect the weights $w_{kj}$ in $w \in \R^{Kd}$. Noting that since $w_{kj} \geq 0$ for all $k \in [K], j \in [d]$,
    \begin{align*}
        \onenorm{w} 
        = \sum_{k=1}^K \sum_{j=1}^d \Abs{w_{kj}}
        = \sum_{k=1}^K \sum_{j=1}^d w_{kj}
        = d \sum_{k=1}^K \sigma_k^2 \eta_k^2 \Pi_{k,K}^{\mu,\nu}. 
    \end{align*}
    Hanson-Wright (Theorem \ref{thm:hansonwright}) and $\infnorm{w} \leq \twonorm{w} \leq \onenorm{w}$ justify the inequalities below:
    \begin{align*}
        \P\Parent{\mathcal E^c}
        &:= \P\Parent{\sum_{k=1}^K \sum_{j=1}^d w_{kj}(X_k)_j^2 \geq \Parent{1+2\sqrt{\log(1/\delta)}+2\log(1/\delta)} \onenorm{w}} \\
        &= \P\Parent{\sum_{k=1}^K \sum_{j=1}^d w_{kj}(X_k)_j^2 \geq \onenorm{w} + 2\onenorm{w} \sqrt{\log(1/\delta)} + 2\onenorm{w} \log(1/\delta)} \\
        &\leq \P\Parent{\sum_{k=1}^K \sum_{j=1}^d w_{kj}((X_k)_j^2-1) \geq 2\twonorm{w} \sqrt{\log(1/\delta)} + 2\infnorm{w} \log(1/\delta)}
        \leq \delta. 
    \end{align*}
    On  $\mathcal E$ we therefore have that 
    \begin{align*}
        T_{\operatorname{DP}}
        &= \sum_{k=1}^K \sum_{j=1}^d \sigma_k^2 \eta_k^2 \Pi_{k,K}^{\mu,\nu} (X_k)_j^2
        \leq \Parent{1+2\sqrt{\log(1/\delta)}+2\log(1/\delta)} d \sum_{k=1}^K \sigma_k^2 \eta_k^2 \Pi_{k,K}^{\mu,\nu} \\
        &\leq 6d \log(1/\delta) \sum_{k=1}^K \sigma_k^2 \eta_k^2 \Pi_{k,K}^{\mu,\nu}. 
    \end{align*}
    To obtain the last inequality we use that $1 \leq \sqrt{2\log(1/\delta)}$ as $\delta < \frac{1}{2} < e^{-1/2}$. 
    
    We will use Theorem \ref{thm:freedman} to control the cross term $T_{\operatorname{Cross}}$. For this, we first verify that it is a sum of martingale differences that are conditionally sub-Gaussian. Second, we show that the sum of the conditional variance proxies can be controlled by a linear combination of the differences themselves. 
    Let $\curly{F_k}_{k=1}^K$ be the filtration where
    \begin{align*}
        \mathcal F_k
        := \sigma\Parent{\theta_0, \curly{X_i : i \in B_1 \cup \dotsm \cup B_k}, Z_1,...,Z_k}.
    \end{align*}
    Corollary \ref{cor:condsubgauss} shows that $M_k$ is a martingale difference sequence \wrt{} $\curly{F_k}_{k=1}^K$ and that 
    \begin{align*}
        \E\Bracket{\exp\Parent{\lambda M_k}|\mathcal F_{k-1}}
        &\leq \exp\Parent{\frac{\lambda^2}{2} \Parent{G^2 + \sigma_k^2} \twonorm{\theta_{k-1}-\theta^\star-\eta_k\Delta_{2,k}}^2}
        \quad \text{for all $\lambda \in \R$}. 
    \end{align*}
    In $T_{\operatorname{Cross}}$ the $M_k$ are rescaled by $2\eta_k\Pi_{k,n}^{\mu,\nu}$. Rescaling maintains the martingale property and the $D_k = 2\eta_k \Pi_{k,n}^{\mu,\nu} M_k$ are still conditionally sub-Gaussian with adjusted variance proxies
    \begin{align*}
        v_{k-1}
        := 4\eta_k^2 \Parent{\Pi_{k,K}^{\mu,\nu}}^2 (G^2+\sigma_k^2) \twonorm{\theta_{k-1}-\theta^\star-\eta_k \Delta_{2,k}}^2. 
    \end{align*}
    The sum of the conditional variance proxies is then given by 
    \begin{align*}
        V_K
        := \sum_{k=1}^K v_{k-1}
        &= \sum_{k=1}^K 4\eta_k^2 \Parent{\Pi_{k,K}^{\mu,\nu}}^2 (G^2+\sigma_k^2) \cdot \twonorm{\theta_{k-1}-\theta^\star-\eta_k \Delta_{2,k}}^2 \\
        &\leq \sum_{k=1}^K 4\eta_k^2 \Parent{\Pi_{k,K}^{\mu,\nu}}^2 (G^2+\sigma_k^2) \cdot \Parent{1-\frac{2\eta_k\mu\nu}{\mu+\nu}}\twonorm{\theta_{k-1}-\theta^\star}^2 \\
        &= \sum_{k=1}^K 4\eta_k^2 \Pi_{k,K}^{\mu,\nu} (G^2+\sigma_k^2) \cdot \Parent{1-\frac{2\eta_k\mu\nu}{\mu+\nu}} \Pi_{k,K}^{\mu,\nu} \twonorm{\theta_{k-1}-\theta^\star}^2. 
    \end{align*}
    Here, the inequality follows from \eqref{eq:deterministic_update}. Corollary \ref{cor:prodsimplification} shows that for any $l \in [k-1]$, 
    \begin{align*}
        \Parent{1-\frac{2\eta_k \mu\nu}{\mu+\nu}} \Pi_{l,k-1}^{\mu,\nu} \Pi_{k,K}^{\mu,\nu}
        = \Pi_{l,K}^{\mu,\nu}. 
    \end{align*}
    We may now apply Lemma \ref{lem:estimationerror} with $K = k-1$ and again use $\twonorm{\Delta_{1,k}} \leq G$ to get
    \begin{align*}
        \twonorm{\theta_{k-1}-\theta^\star}^2 
        &\leq \twonorm{\theta_0-\theta^\star}^2 \cdot \Pi_{0,k-1}^{\mu,\nu} + \sum_{j=1}^{k-1} 2\eta_j^2\Parent{G^2 + \twonorm{Z_j}^2} \cdot \Pi_{j,k-1}^{\mu,\nu} + \sum_{j=1}^{k-1} 2\eta_j M_j \cdot \Pi_{j,k-1}^{\mu,\nu}.
    \end{align*}
    Combining both of the results above and using that the summands are non-negative, 
    \begin{align*}
        \Parent{1-\frac{2\eta_k\mu\nu}{\mu+\nu}} & \Pi_{k,K}^{\mu,\nu} \twonorm{\theta_{k-1}-\theta^\star}^2 \\
        &\leq \twonorm{\theta_0-\theta^\star}^2 \cdot \Pi_{0,K}^{\mu,\nu} + \sum_{j=1}^{k-1} 2\eta_j^2\Parent{G^2 + \twonorm{Z_j}^2} \cdot \Pi_{j,K}^{\mu,\nu} + \sum_{j=1}^{k-1} 2\eta_j M_j \cdot \Pi_{j,K}^{\mu,\nu} \\
        &\leq T_{\operatorname{Bias}} + T_{\operatorname{SGD}} + T_{\operatorname{DP}} + \sum_{j=1}^{k-1} 2\eta_j M_j \cdot \Pi_{j,K}^{\mu,\nu}. 
    \end{align*}
    Note here, that $T_{\operatorname{Bias}}$, $T_{\operatorname{SGD}}$ and $T_{\operatorname{DP}}$ do not depend on the index $k-1$. Plugging the expression back into the upper bound on $V_K$, we get
    \begin{align*}
        V_K
        &\leq \sum_{k=1}^K 4\eta_k^2 \Pi_{k,K}^{\mu,\nu} (G^2+\sigma_k^2) \cdot \Parent{1-\frac{2\eta_k\mu\nu}{\mu+\nu}} \Pi_{k,K}^{\mu,\nu} \twonorm{\theta_{k-1}-\theta^\star}^2 \\
        &\leq \sum_{k=1}^K 4\eta_k^2 \Pi_{k,K}^{\mu,\nu} (G^2+\sigma_k^2) \cdot \Parent{T_{\operatorname{Bias}} + T_{\operatorname{SGD}} + T_{\operatorname{DP}} + \sum_{j=1}^{k-1} 2\eta_j M_j \Pi_{j,K}^{\mu,\nu}} \\
        &= \sum_{k=1}^K 4\eta_k^2 \Pi_{k,K}^{\mu,\nu} (G^2+\sigma_k^2) \sum_{j=1}^{k-1} 2\eta_j M_j \Pi_{j,K}^{\mu,\nu} + \beta(\delta), \\
        &\text{where $\beta(\delta)
        := \Parent{T_{\operatorname{Bias}} + T_{\operatorname{SGD}} + T_{\operatorname{DP}}} \sum_{k=1}^K 4\eta_k^2 \Pi_{k,K}^{\mu,\nu} (G^2+\sigma_k^2)$}. 
    \end{align*}
    By changing the order of summation, we rewrite this into a bound as required by the Generalized Freedman's inequality (Theorem \ref{thm:freedman}):  
    \begin{align*}
        V_K
        &\leq \sum_{k=1}^K 4\eta_k^2 \Pi_{k,K}^{\mu,\nu} (G^2+\sigma_k^2) \sum_{j=1}^{k-1} 2\eta_j M_j \cdot \Pi_{j,K}^{\mu,\nu} + \beta(\delta) \\
        &= \sum_{j=1}^K \underbrace{\Parent{\sum_{k=j+1}^K 4 \eta_k^2 \Pi_{k,K}^{\mu,\nu} (G^2+\sigma_k^2)}}_{=: \alpha_j} \cdot 2\eta_j \Pi_{j,K}^{\mu,\nu} M_j + \beta(\delta).
    \end{align*}
    We proceed by bounding $\infnorm{\alpha}$ and $\beta(\delta)$. As the summands in $\alpha_j$ are non-negative,
    \begin{align*}
        \infnorm{\alpha}
        &
        = \max_{j \in [K]} \alpha_j
        = \max_{j \in [K]} \underbrace{\sum_{k=j+1}^K 4 \eta_k^2 \Pi_{k,K}^{\mu,\nu} (G^2+\sigma_k^2)}_{\alpha_j}
        \leq \sum_{k=1}^K 4 \eta_k^2 \Pi_{k,K}^{\mu,\nu} (G^2+\sigma_k^2)
        =: \infnorm{\bar\alpha}.  
    \end{align*}
    When expressing $\beta(\delta)$ in terms of $\infnorm{\bar\alpha}$, on $\mathcal E$ the following inequality holds: 
    \begin{align*}
        \beta(\delta) 
        &= \Parent{T_{\operatorname{Bias}} + T_{\operatorname{SGD}} + T_{\operatorname{DP}}} \cdot \infnorm{\bar\alpha}
        \leq \frac{1}{2} \Parent{T_{\operatorname{Bias}} + T_{\operatorname{SGD}} + T_{\operatorname{DP}}}^2 + \frac{\infnorm{\bar\alpha}^2}{2} \\
        &\leq \frac{1}{2} \Parent{\twonorm{\theta_0-\theta^\star}^2 \Pi_{0,K}^{\mu,\nu} + 2G^2 \sum_{k=1}^K \eta_k^2 \Pi_{k,K}^{\mu,\nu} + 6d\log(1/\delta) \sum_{k=1}^K \sigma_k^2 \eta_k^2 \Pi_{k,K}^{\mu,\nu}}^2 + \frac{\infnorm{\bar\alpha}^2}{2}. 
    \end{align*}
    Here, we use Young's inequality and plug in $T_{\operatorname{Bias}}$ and the upper bounds on $T_{\operatorname{SGD}}$ and $T_{\operatorname{DP}}$. 
    
    Using the sub-additivity of $x \mapsto \sqrt{x}$, 
    \begin{align*}
        \sqrt{2\bar\beta(\delta)}
        &\leq \twonorm{\theta_0-\theta^\star}^2 \Pi_{0,K}^{\mu,\nu} + 2G^2 \sum_{k=1}^K \eta_k^2 \Pi_{k,K}^{\mu,\nu} + 6d\log(1/\delta) \sum_{k=1}^K \sigma_k^2 \eta_k^2 \Pi_{k,K}^{\mu,\nu} + \infnorm{\bar\alpha} \\
        &=: \sqrt{2\bar\beta(\delta)}
    \end{align*}

    We have shown that on $\mathcal E$ and thus with probability at least $1-\delta$, $V_K$ is upper bounded by a linear combination of the martingales $2\eta_j \Pi_{j,n}^{\mu,\nu} M_j$ with weights $\alpha_j \leq \infnorm{\bar\alpha}$ and bias $\bar\beta(\delta)$. We may apply Freedman's inequality (Theorem \ref{thm:freedman}). With probability at least $1-2\delta$, 
    \begin{align*}
        T_{\operatorname{Cross}} 
        = \sum_{k=1}^K 2\eta_k M_k \cdot \Pi_{k,K}^{\mu,\nu}
        &\leq \infnorm{\bar\alpha} \cdot 4\log(1/\delta) + \sqrt{\bar\beta(\delta) \cdot 8\log(1/\delta)} \\
        &\leq \Parent{\infnorm{\bar\alpha} + \sqrt{\bar\beta(\delta)}} 4\log(1/\delta).
    \end{align*}
    Here, for the second inequality we use $\delta < \frac{1}{2} < e^{-1/2}$. Collecting terms, 
    \begin{align*}
        &\twonorm{\theta_n-\theta^\star}^2
        \leq T_{\operatorname{Bias}} + T_{\operatorname{SGD}} + T_{\operatorname{DP}} + T_{\operatorname{Cross}} \\
        &\quad \lesssim \Parent{\twonorm{\theta_0-\theta^\star}^2 \Pi_{0,K}^{\mu,\nu} + G^2 \sum_{k=1}^K \eta_k^2 \Pi_{k,K}^{\mu,\nu} +  d\log(1/\delta) \sum_{k=1}^K \sigma_k^2 \eta_k^2 \Pi_{k,K}^{\mu,\nu} + \infnorm{\bar \alpha}} \log(1/\delta). 
    \end{align*}
    The statement follows as $\infnorm{\bar\alpha}$ is dominated in order by the second and third term above. 
\end{proof}

\begin{corollary}
    \label{cor:condsubgauss}
    Let $\curly{\theta_k}_{k=1}^K$ be the iterates in Algorithm \ref{alg:noisySGD} and $\curly{F_k}_{k=1}^K$ be the filtration with
    \begin{align*}
        \mathcal F_k
        := \sigma\Parent{\theta_0, \curly{X_i : i \in B_1 \cup \dotsm \cup B_k}, Z_1,...,Z_k}.
    \end{align*}
    Then, $M_k := -(\theta_{k-1}-\theta^\star-\eta_k \Delta_{2,k})\T(\Delta_{1,k}+Z_k)$ is a martingale difference sequence \wrt{} the filtration above. If we further suppose that $\twonorm{\Delta_{1,k}} \leq G$, it holds that
    \begin{align*}
        \E\Bracket{\exp\Parent{\lambda M_k}|\mathcal F_{k-1}}
        &\leq \exp\Parent{\frac{\lambda^2}{2} \Parent{G^2 + \sigma_k^2} \twonorm{\theta_{k-1}-\theta^\star-\eta_k\Delta_{2,k}}^2}
        \quad \text{for all $\lambda \in \R$}. 
    \end{align*}
\end{corollary}

\begin{proof}
    For convenience, define the following quantity 
    \begin{align*}
        \delta_k 
        &:= \theta_{k-1}-\theta^\star-\eta_k \Delta_{2,k} \\
        &= \theta_{k-1}-\theta^\star-\eta_k \Parent{\nabla_\theta \mathcal L_P(\theta_{k-1}) - \nabla_\theta \mathcal L_P(\theta^\star)}. 
    \end{align*}
    We first show that $\curly{M_k}_{k=1}^K$ is a martingale difference sequence. For all $k \in [K]$, 
    \begin{align*}
        \E\Bracket{M_k|\mathcal F_{k-1}}
        &= \E\Bracket{\delta_k\T (\Delta_{1,k} + Z_k)|\mathcal F_{k-1}} \\
        &= \delta_k\T \Parent{\E\Bracket{\Delta_{1,k}|\mathcal F_{k-1}} + \E\Bracket{Z_k|\mathcal F_{k-1}}}
        = 0.
    \end{align*}
    Here, $\E\Bracket{Z_k|\mathcal F_{k-1}} = 0$ as $Z_k \sim \mathcal N(0,\sigma_k^2 I_d)$ is independent of $\theta_0,...,\theta_{k-1}$. As the $X_i$ are identically distributed, the batches $\curly{B_k}_{k=1}^K$ are disjoint and $\theta_{k-1}$ and $\curly{X_i}_{i \in B_k}$ are mutually independent, Corollary \ref{cor:iidbatchriskisunbiased} applied conditionally yields that $\E\Bracket{\Delta_{1,k}|\mathcal F_{k-1}} = 0$.

    We proceed by showing the conditional sub-Gaussianity of $M_k$ given $\mathcal F_{k-1}$. For all $\lambda \in \R$, 
    \begin{align*}
        \E\Bracket{\exp\Parent{\lambda M_k}|\mathcal F_{k-1}}
        &= \E\Bracket{\exp\Parent{-\lambda \parent{\theta_{k-1}-\theta^\star-\eta_k \Delta_{2,k}}\T(\Delta_{1,k} + Z_k)}\Big|\mathcal F_{k-1}} \\
        &= \E\Bracket{\exp\Parent{-\lambda \delta_k \T(\Delta_{1,k} + Z_k)}|\mathcal F_{k-1}} \\
        &\overset{(i)}= \E\Bracket{\exp\Parent{-\lambda \delta_k\T \Delta_{1,k}}|\mathcal F_{k-1}} \E\Bracket{\exp\Parent{-\lambda \delta_k\T Z_k}|\mathcal F_{k-1}} \\
        &\overset{(ii)}\leq \exp\Parent{\frac{\lambda^2}{2} G^2 \twonorm{\delta_k}^2} \exp\Parent{\frac{\lambda^2}{2} \sigma_k^2 \twonorm{\delta_k}^2} \\
        &\leq \exp\Parent{\frac{\lambda^2}{2} \Parent{G^2 + \sigma_k^2} \twonorm{\delta_k}^2}. 
    \end{align*}
    Above, $(i)$ holds by the independence of $B_k, \theta_{k-1}$ and $Z_k$. It remains to justify the inequality $(ii)$. As we assume that the stochastic gradients have norm bounded by $G$, 
    \begin{align*}
        \Abs{\delta_k\T \Delta_{1,k}}
        \leq \twonorm{\delta_k} \twonorm{\Delta_{1,k}}
        \leq G \twonorm{\delta_k}. 
    \end{align*}
    Thus, by Hoeffding's Lemma applied conditionally, $\delta_k\T \Delta_{1,k}$ is $G^2\twonorm{\delta_k}^2$-sub-Gaussian. Since $Z_k \sim \mathcal N(0, \sigma_k^2)$ it holds that $\delta_k\T Z_k \sim \mathcal N(0, \sigma_k^2 \twonorm{\delta_k}^2)$ given $\mathcal F_{k-1}$. This holds by the affine transformation theorem and because Gaussians are fully determined by their mean and variance, which here are 
    \begin{align*}
        \E\Bracket{\delta_k\T Z_k |\mathcal F_{k-1}}
        &= \delta_k\T \E\Bracket{Z_k|\mathcal F_{k-1}}
        = 0, \\
        \var{\delta_k\T Z_k|\mathcal F_{k-1}}
        &= \delta_k\T \var{Z_k|\mathcal F_{k-1}} \delta_k
        = \delta_k\T \sigma_k^2 I_d \delta_k
        = \sigma_k^2 \twonorm{\delta_k}^2. 
    \end{align*}
    Accordingly, $\delta_k\T Z_k$ is $\sigma_k^2 \twonorm{\delta_k}^2$-sub-Gaussian. The inequality $(ii)$ and thus also the statement are now recovered using the definition of sub-Gaussianity in terms of an MGF bound. 
\end{proof}

\subsubsection{Proof of Theorem \ref{thm:esterrorinprob}}

\begin{proof}
    We first verify that $\eta_k \leq \frac{2}{\mu+\nu_\ell} \leq \frac{2}{\nu_\ell}$ for all $k \in [n]$. For this need to check that 
    \begin{align*}
        \max_{k \in [n]} \eta_k 
        = \eta_1 
        = \frac{\mu+\nu_\ell}{\mu\nu_\ell} \frac{2}{1+\lceil\nu_\ell/\mu+2\rceil}
       \leq \frac{2}{\mu+\nu_\ell}
        \quad \Leftrightarrow \quad 
        (\mu+\nu_\ell)^2
        \leq \mu\nu_\ell + \mu\nu_\ell \lceil\nu_\ell/\mu+2\rceil. 
    \end{align*}
    The latter inequality is implied by the following statements: 
    \begin{align*}
      \mu
        \leq \nu_\ell 
        \quad \Leftrightarrow \quad 
        \mu^2+\nu_\ell^2+2\mu\nu_\ell  
       \leq \mu\nu_\ell + \nu_\ell^2 + 2\mu\nu_\ell
        \quad \Leftrightarrow 
        (\mu+\nu_\ell)^2 
      \leq \mu\nu_\ell + \mu\nu_\ell (\nu_\ell/\mu+2)
    \end{align*}
    Algorithm \ref{alg:noisySGD} with $b = 1$ is private by the reasoning in Theorem \ref{thm:esterrorinexp}.

    Moreover, due to the uniform bound $\sup_{\theta \in \Theta, x_i \in \mathcal X} \twonorm{\nabla_\theta \ell(\theta, x_i)} \leq \frac{L}{2}$ it holds that 
    \begin{align*}
        \twonorm{\Delta_{1,k}}^2
        = \twonorm{\nabla_\theta \mathcal L_P(\theta_{k-1}) - \nabla_\theta \mathcal L_{B_k}(\theta_{k-1})}^2
        \leq 2\twonorm{\nabla_\theta \mathcal L_P(\theta_{k-1})}^2 + 2\twonorm{\nabla_\theta \mathcal L_{B_k}(\theta_{k-1})}^2
        \leq L^2 
        =: G^2. 
    \end{align*}
    Since $\eta_k \leq \frac{2}{\mu+\nu_\ell}$, by Lemma \ref{lem:esterrorinprobgen} with $G = L$ we get
    \begin{align*}
        \twonorm{\theta_K-\theta^\star}^2
        &\lesssim \Parent{\twonorm{\theta_0-\theta^\star}^2 \cdot \Pi_{0,K}^{\mu,\nu_\ell} + L^2 \sum_{k=1}^K \eta_k^2 \Pi_{k,K}^{\mu,\nu_\ell} + d\log(1/\delta) \sum_{k=1}^K \sigma_k^2 \eta_k^2 \Pi_{k,K}^{\mu,\nu_\ell}} \log(1/\delta). 
    \end{align*}
    We now control the three terms above in order. 
    
    For the bias term $T_{\operatorname{Bias}}$, writing $\beta:=4$ and $\tau := \lceil \nu_\ell/\mu+2\rceil$, Corollary \ref{cor:prodBias} and $\frac{\tau}{n}\leq 1$ yield 
    \begin{align*}
        T_{\operatorname{Bias}}
        := \twonorm{\theta_0-\theta^\star}^2 \cdot \Pi_{0,n}^{\mu,\nu_\ell}
        \leq \twonorm{\theta_0-\theta^\star}^2 \cdot \Parent{\frac{\tau}{n}}^\beta
        \leq \twonorm{\theta_0-\theta^\star}^2 \cdot \Parent{\frac{\tau}{n}}
        \lesssim \twonorm{\theta_0-\theta^\star}^2 \cdot \frac{\nu_\ell}{\mu n}
    \end{align*}
    By Corollary \ref{cor:sumSGD} and using that $\frac{\mu+\nu_\ell}{2\mu\nu_\ell} \leq \frac{1}{\mu}$ we further have
    \begin{align*}
        T_{\operatorname{SGD}}
        := L^2 \sum_{k=1}^n \eta_k^2 \Pi_{k,n}^{\mu,\nu_\ell}
        \leq \frac{L^2 \beta^2}{n} \Parent{\frac{\mu+\nu}{2\mu\nu}}^2
        \leq \frac{L^2\beta^2}{\mu^2 n}.
    \end{align*}
    Finally, Corollary \ref{cor:sumDP} lets us upper bound the privacy noise term $T_{\operatorname{DP}}$: 
    \begin{align*}
        T_{\operatorname{DP}}
        := d\log(1/\delta) \sum_{k=1}^n \sigma_k^2 \eta_k^2 \Pi_{k,n}^{\mu,\nu_\ell}
        \leq d\log(1/\delta) \frac{3\alpha L^2\beta^2}{\epsilon n^2} \Parent{\frac{\mu+\nu_\ell}{2\mu\nu_\ell}}^2
        \leq \frac{3\alpha dL^2\beta^2}{\mu^2 \epsilon n^2} \log(1/\delta). 
    \end{align*}
    Collecting terms recovers the statement.
\end{proof}

\section{Proofs of Section \ref{sec:minimax}}
\label{app:minimax}

\begin{lemma}{(Proposition 3.2, \cite{cai:wang:zhang2025}).}
    \label{lem:lbexpscoreattack}
    Suppose $(Y, X) \in \mathcal Z^n \subseteq \R^{n \times (d+1)}$ contains \iid{} samples $(Y_i, X_i)\in \mathcal Z \subseteq \R^{d+1}$ from $\operatorname{GLM}_{\Theta_0}^\psi(\vartheta)$. Let $V := (V_1,...,V_d)\T \in [-1,1]^d$ have \iid{} entries $V_j$ with density $p(v) = 15(1-v^2)^2/16$. The distribution $\pi$ of $\vartheta = a V$ is supported on $\Theta_0$ if $a := R_\Theta/(2\sqrt{d})$. Then, given that $M \in \mathcal M$ fulfills $\E\Bracket{\E\Bracket{\twonorm{M((Y,X))-\vartheta}^2|\vartheta}} \leq \frac{R_\Theta^2}{160}$,
    \begin{align*}
        \frac{d}{2}
        \leq \sum_{i=1}^n \E\Bracket{\E\Bracket{\mathcal A_\vartheta\Parent{(Y_i,X_i), M((Y,X))}\Big|\vartheta}}. 
    \end{align*}
\end{lemma}

\subsection{Proofs for Score Attack}

\subsubsection{Proof of Lemma \ref{lem:ubexpscoreattack}}

\begin{proof}
    Define the family of distributions $\mathcal Q := \curly{Q_z : z \in \mathcal Z^n}$ s.t.\ $M(z) \sim Q_z$ for all $z \in \mathcal Z^n$, meaning that $M(Z) \sim Q_Z$ and $M(Z_{(i)}) \sim Q_{Z_{(i)}}$. Now, define the likelihood ratio 
    \begin{align*}
        L_i(m)
        := \frac{\mathrm{d}Q_Z(m)}{\mathrm{d}Q_{Z_{(i)}}}
        \quad \text{for all $m \in \R^d$}.
    \end{align*}
    For all $\tilde z \in \mathcal Z$ we can then perform the following change of measure:
    \begin{align*}
        \E\bracket{\mathcal A_\theta\parent{\tilde z,M(Z)}|Z,Z_i^\prime}
        = \int_{\R^d} \mathcal A_\theta\parent{\tilde z,m} \mathrm{d}Q_Z(m) 
        = \int_{\R^d} \mathcal A_\theta\parent{\tilde z,m} L_i(m) \mathrm{d}Q_{Z_{(i)}}(m).
    \end{align*}
    Using this identity and the Cauchy Schwarz inequality gets us 
    \begin{align*}
        \E\Bracket{\mathcal A_\theta\Parent{Z_i,M(Z)}|Z,Z_i^\prime} &- \E\Bracket{\mathcal A_\theta\Parent{Z_i,M(Z_{(i)})}|Z,Z_i^\prime} \\
        &= \int_{\R^d} \Parent{L_i(m)-1}\mathcal A_\theta(Z_i,m) \mathrm{d} Q_{Z_{(i)}}(m) \\
        &\leq \sqrt{\int_{\R^d} (L_i(m)-1)^2 \mathrm{d}Q_{Z_{(i)}}(m)} \sqrt{\int_{\R^d} \mathcal A_\theta(Z_i,m)^2 \mathrm{d}Q_{Z_{(i)}}(m)} \\
        &= \sqrt{\chi^2\parent{Q_Z||Q_{Z_{(i)}}}} \sqrt{\E\Bracket{\mathcal A_\theta\Parent{Z_i,M(Z_{(i)})}^2|Z, Z_i^\prime}} \\
        &\leq C_{\alpha, \epsilon} \sqrt{\E\Bracket{\mathcal A_\theta\Parent{Z_i,M(Z_{(i)})}^2|Z, Z_i^\prime}}.
    \end{align*}
    The last inequality uses $\chi^2\parent{Q_Z||Q_{Z_{(i)}}} \leq e^\epsilon-1 =: C_{\alpha,\epsilon}^2$ for all $\alpha \geq 2$ which holds because $D_2(P||Q)=\log(1+\chi^2\parent{P||Q})$ and $D_2(P||Q)\leq D_\alpha(P||Q)\leq \epsilon$ for $\alpha\geq 2$. Applying both statements in Theorem \ref{thm:soundness} and Jensen's inequality to the concave function $x \mapsto \sqrt{x}$,
    \begin{align*}
        \E\Bracket{\mathcal A_\theta\Parent{Z_i,M(Z)}}
        &= \E\Bracket{\mathcal A_\theta\Parent{Z_i,M(Z)}} - \E\Bracket{\mathcal A_\theta\Parent{Z_i,M(Z_{(i)})}} \\
        &\leq C_{\alpha,\epsilon} \E\Bracket{\sqrt{\E\Bracket{\mathcal A_\theta\Parent{Z_i,M(Z_{(i)})}^2|Z, Z_i^\prime}}} \\
        &\leq C_{\alpha,\epsilon} \sqrt{\E\Bracket{\mathcal A_\theta\Parent{Z_i,M(Z_{(i)})}^2}} \\
        &\leq C_{\alpha,\epsilon} \cdot \sqrt{\eigmax\Parent{I(\theta)} \cdot \E\Bracket{\twonorm{M(Z)-\theta}^2}}. 
    \end{align*}
    Plugging in the definition of $C_{\alpha,\epsilon}$ recovers the statement. 
\end{proof}

\subsection{Statistical Lower Bound}

\subsubsection{Proof of Lemma \ref{lem:lbpriv}}

\begin{proof}
    Fix an estimator $M \in \mathcal M_{\alpha,\epsilon}$ with $\mathcal Z := \mathcal X \times \mathcal Y \subseteq \R^{d+1}$ and $\Theta_0 \subset \R^d$. When 
    \begin{align*}
        \sup_{\vartheta \in \Theta_0} \E\Bracket{\twonorm{M((Y,X))-\vartheta}^2|\vartheta} > R_\Theta^2/160
    \end{align*}
    the assumption that $\frac{d^2}{\epsilon n^2} \leq 1$ yields
    \begin{align*}
        \sup_{\vartheta \in \Theta_0} \E\Bracket{\twonorm{M((Y,X))-\vartheta}^2|\vartheta}
        > \frac{R_\Theta^2}{160}
        \geq \frac{R_\Theta^2}{160} \frac{d^2}{\epsilon n^2}
        \geq \frac{R_\Theta^2}{160} \frac{d^2}{(e^\epsilon-1) n^2}.
    \end{align*}
    On the complement event, by definition of the supremum, we instead have 
    \begin{align*}
        \E\Bracket{\E\Bracket{\twonorm{M((Y,X))-\vartheta}^2|\vartheta}}
        \leq \sup_{\vartheta \in \Theta_0} \E\Bracket{\twonorm{M((Y,X))-\vartheta}^2|\vartheta}
        \leq \frac{R_\Theta^2}{160}. 
    \end{align*}
    As the GLM is in canonical form and hence $\E\Bracket{Y_i|X_i, \vartheta} = \psi^\prime(X_i\T\vartheta)$ we have $\E\Bracket{S_\vartheta(Z_i)|\vartheta} = 0$. Using that $\infnorm{\psi^{\prime\prime}} \leq C_\psi^{\prime\prime}$ and $\E\Bracket{X_i X_i\T} \preceq \xi I_d$, we can control the Fisher information as 
    \begin{equation}
        I(\vartheta) 
        = - \E\Bracket{\nabla_\theta^2 \log\Parent{p_\vartheta(Y_1,X_1)}}
        = \E\Bracket{\frac{\psi^{\prime\prime}(X_1\T\vartheta)}{c(\varsigma)} X_1 X_1\T}
        \preceq \frac{C_\psi^{\prime\prime}}{c(\varsigma)} \E\Bracket{X_1 X_1\T}
        \preceq \frac{C_\psi^{\prime\prime}\xi}{c(\varsigma)} I_d. 
        \label{eq:bound_Fisher_info}
    \end{equation}
    Thus, the conditions of Lemma \ref{lem:ubexpscoreattack} are fulfilled. By Lemma \ref{lem:lbexpscoreattack} and \ref{lem:ubexpscoreattack}, for all $\vartheta \in \Theta_0$
    \begin{align*}
        \frac{d}{2}
        &\leq \sum_{i=1}^n \E\Bracket{\E\Bracket{\mathcal A_\vartheta\Parent{(Y_i,X_i), M((Y,X))}\Big|\vartheta}} \\
        &\leq \sqrt{e^\epsilon-1} \cdot \sum_{i=1}^n \E\Bracket{\sqrt{\eigmax\Parent{I(\vartheta))} \cdot \E\Bracket{\twonorm{M((Y,X))-\vartheta}^2|\vartheta}}} \\
        &\leq \sqrt{e^\epsilon-1} \cdot \sum_{i=1}^n \E\Bracket{\sqrt{\frac{C_\psi^{\prime\prime} \xi}{c(\varsigma)} \cdot \E\Bracket{\twonorm{M((Y,X))-\vartheta}^2|\vartheta}}} \\
        &= n \sqrt{e^\epsilon-1} \sqrt{\frac{C_\psi^{\prime\prime} \xi}{c(\varsigma)}} \cdot \E\Bracket{\sqrt{\E\Bracket{\twonorm{M((Y,X))-\vartheta}^2|\vartheta}}} \\
        &\leq n \sqrt{e^\epsilon-1} \sqrt{\frac{C_\psi^{\prime\prime} \xi}{c(\varsigma)}} \cdot \sqrt{\E\Bracket{\E\Bracket{\twonorm{M((Y,X))-\vartheta}^2|\vartheta}}}. 
    \end{align*}
    Here, the last inequality is concave Jensen's inequality applied to $x \mapsto \sqrt{x}$. Solving for the expectation and again by definition of the supremum, we then have
    \begin{align*}
        \sup_{\vartheta \in \Theta_0} \E\Bracket{\twonorm{M((Y,X))-\vartheta}^2|\vartheta} 
        \geq \E\Bracket{\E\Bracket{\twonorm{M((Y,X))-\vartheta}^2|\vartheta}}
        \geq \frac{c(\varsigma)}{4 C_\psi^{\prime\prime} \xi} \frac{d^2}{(e^\epsilon-1) n^2}
    \end{align*}
    As either of the lower bounds hold for all $M \in \mathcal M_{\alpha,\epsilon}$, we recover the statement. 
\end{proof}

\subsubsection{Proof of Theorem \ref{thm:lbcombined}}

\begin{proof}
    As $P_{\mathcal X}$ and its density $p(x)$ do not depend on $\vartheta$, 
    \begin{align*}
        S_\vartheta((y,x))
        &= \nabla_\vartheta \log\Parent{p_\vartheta(y,x)}
        = \nabla_\vartheta \log\Parent{p_\vartheta(y|x)} \\
        &= \nabla_\vartheta \Parent{\frac{y x\T\vartheta-\psi(x\T\vartheta)}{c(\varsigma)}}
        = \frac{\Parent{y-\psi^\prime(x\T\vartheta)}}{c(\varsigma)} \cdot x. 
    \end{align*}
    For $\epsilon \lesssim 1$ we have $e^\epsilon-1 \lesssim \epsilon$. Combining this with the bound on the maximum eigenvalue of $I(\vartheta)$ in \eqref{eq:bound_Fisher_info}, the private lower bound in Lemma \ref{lem:lbpriv} simplifies into
    \begin{align*}
        \inf_{M \in \mathcal M_{\alpha,\epsilon}}\sup_{\vartheta \in \Theta_0} \E\Bracket{\twonorm{M((Y,X))-\vartheta}^2|\vartheta} 
        &\gtrsim \Parent{R_\Theta^2 \wedge \frac{1}{\eigmax(I(\theta))}} \cdot \frac{d^2}{(e^\epsilon-1) n^2} \\
        &\gtrsim \Parent{R_\Theta^2 \wedge \frac{c(\varsigma)}{C_\psi^{\prime\prime} \xi}} \cdot \frac{d^2}{(e^\epsilon-1) n^2} \\
        &\gtrsim \Parent{R_\Theta^2 \wedge \frac{c(\varsigma)}{C_\psi^{\prime\prime} \xi}} \cdot \frac{d^2}{\epsilon n^2}. 
    \end{align*}
    At the same time, the non-private lower bound in Lemma \ref{lem:lbnonpriv} is applicable as $\frac{d}{n} \leq \sqrt{\epsilon} \leq 1$. Since the Renyi-DP estimators are a subset of all estimators, we have $\mathcal M_{\alpha,\epsilon} \subset \mathcal M$ and
    \begin{align*}
        \inf_{M \in \mathcal M_{\alpha,\epsilon}}\sup_{\vartheta \in \Theta_0} \E\Bracket{\twonorm{M((Y,X))-\vartheta}^2|\vartheta} 
        &\geq \inf_{M \in \mathcal M}\sup_{\vartheta \in \Theta_0} \E\Bracket{\twonorm{M((Y,X))-\vartheta}^2|\vartheta} \\
        &\gtrsim \Parent{R_\Theta^2 \wedge \frac{c(\varsigma)}{C_\psi^{\prime\prime} \xi}} \cdot \frac{d}{n}. 
    \end{align*}
    Combining both bounds yields the statement. 
\end{proof}

\begin{lemma}
    \label{lem:lbnonpriv}
    Let $(Y, X) \in \mathcal Z^n \subseteq \R^{n \times (d+1)}$ contain \iid{} samples $(Y_i, X_i\T)\T \in \mathcal Z \subseteq \R^{d+1}$ from $\operatorname{GLM}_{\Theta_0}^\psi(\vartheta)$ where $\vartheta$ is as in the supremum. Assuming that $\frac{d}{n} \leq 1$, it holds that
    \begin{align*}
        \inf_{M \in \mathcal M} \sup_{\vartheta \in \Theta_0} \E\Bracket{\twonorm{M((Y,X))-\vartheta}^2|\vartheta}
        \gtrsim \Parent{R_\Theta^2 \wedge \frac{c(\varsigma)}{C_\psi^{\prime\prime} \xi}} \cdot \frac{d}{n}. 
    \end{align*}
\end{lemma}

\begin{proof}
    See the proof of \cite[Theorem 3.1]{cai:wang:zhang2025} or \cite[Section 8]{duchi:2023}.
\end{proof}

\subsection{Statistical Upper Bound}

\begin{lemma}
    \label{lem:recoveryPRM}
    Let $(Y, X) \in \mathcal Z^n \subseteq \R^{n \times (d+1)}$ contain \iid{} samples $(Y_i, X_i\T)\T \in \mathcal Z \subseteq \R^{d+1}$ from $\operatorname{GLM}_{\Theta_0}^\psi(\vartheta)$. Then, $\vartheta$ is a minimizer of the negative log-likelihood loss, i.e., 
    \begin{align*}
        \vartheta 
        \in \arg\min_{\theta \in \Theta_0} \E\Bracket{-\log\Parent{p_\theta(Y_1|X_1)}}.
    \end{align*}
\end{lemma}

\subsubsection{Proof of Lemma \ref{lem:esterrorinexpGLMs}}

\begin{proof}
    For the given GLM, the negative log-likelihood loss component is
    \begin{align*}
        \ell(\theta, z_i)
        = -\log\Parent{p_\theta(y_i|x_i)}
        = \frac{\psi(x_i\T \theta) - y_i x_i\T\theta}{c(\varsigma)} - \log(h(y_i,\varsigma)).
    \end{align*}
    Accordingly, the gradient and Hessian can be computed as
    \begin{align*}
        \nabla_\theta \ell(\theta, z_i)
        &= \frac{1}{c(\varsigma)} \Parent{\psi^\prime(x_i\T \theta) - y_i} \cdot x_i
        \quad \text{and} \quad 
        \nabla_\theta^2 \ell(\theta, z_i)
        = \frac{1}{c(\varsigma)} \psi^{\prime\prime}(x_i\T\theta) \cdot x_i x_i\T. 
    \end{align*}
    Given that $\infnorm{\psi^\prime} \leq C_\psi^\prime$ and $\twonorm{X_i} \leq \frac{R_{\mathcal X}}{2} \sqrt{d}$ for all $i \in [n]$, it holds that
    \begin{align*}
        \sup_{\theta \in \Theta_0, Z_i \in \mathcal Z} \twonorm{\nabla_\theta \ell(\theta, Z_i)}
        &= \sup_{\theta \in \Theta_0, Z_i \in \mathcal Z} \frac{1}{c(\varsigma)} \Twonorm{\Parent{\psi^\prime(X_i\T \theta)-Y_i} \cdot X_i} \\
        &\leq \frac{1}{c(\varsigma)} \sup_{\theta \in \Theta_0, Z_i \in \mathcal Z} \Parent{C_\psi^\prime + \abs{Y_i}} \twonorm{X_i} \\
        &\leq \frac{1}{c(\varsigma)}\Parent{C_\psi^\prime + R_{\mathcal Y}} \sup_{\theta \in \Theta_0, Z_i \in \mathcal Z} \twonorm{X_i} \\
        &\leq \Parent{C_\psi^\prime + R_{\mathcal Y}} \frac{R_{\mathcal X}}{2c(\varsigma)} \sqrt{d}.
    \end{align*}
    Thus, the gradients are uniformly bounded with constant $L := \parent{C_\psi^\prime + R_{\mathcal Y}} \frac{R_{\mathcal X}}{c(\varsigma)} \sqrt{d}$.

    Moreover, because $\infnorm{\psi^{\prime\prime}} \leq C_\psi^{\prime\prime}$ and $\twonorm{X_i} \leq \frac{R_{\mathcal X}}{2} \sqrt{d}$ for all $i \in [n]$, we get that
    \begin{align*}
        \nabla_\theta^2 \ell(\theta, Z_i)
        = \frac{1}{c(\varsigma)} \psi^{\prime\prime}(X_i\T\theta) \cdot X_i X_i\T
        \preceq \frac{1}{c(\varsigma)} C_\psi^{\prime\prime} \twonorm{X_i}^2 I_d
        \preceq \frac{1}{4c(\varsigma)} C_\psi^{\prime\prime}R_{\mathcal X}^2d \cdot I_d. 
    \end{align*}
    Therefore, the loss components are also $\nu_\ell$-smooth with $\nu_\ell := \frac{1}{4c(\varsigma)}C_\psi^{\prime\prime} R_{\mathcal X}^2d$. 

    To apply Theorem \ref{thm:esterrorinexp} we now show that the population risk $\mathcal L_P$ is $\mu$-strongly convex with $\mu$ independent of $n,d$. The Hessian of $\mathcal L_n$ corresponding to the design $X \in \R^{n \times d}$ is 
    \begin{align*}
        \nabla_\theta^2 \mathcal L_n(\theta) 
        = \frac 1 n \sum_{i=1}^n \nabla_\theta^2 \ell(\theta, Z_i)
        = \frac 1 n \sum_{i=1}^n \frac{1}{c(\varsigma)} \psi^{\prime\prime}(X_i\T \theta) \cdot X_i X_i\T. 
    \end{align*}
    From this, we compute the Hessian of the population risk $\mathcal L_P$ as 
    \begin{align*}
        \nabla_\theta^2 \mathcal L_P(\theta) 
        = \E\Bracket{\nabla_\theta^2 \mathcal L_n(\theta)}
        = \frac{1}{n} \sum_{i=1}^n \E\Bracket{\frac{1}{c(\varsigma)} \psi^{\prime\prime}(X_i\T\theta) \cdot X_iX_i\T}
        = \frac{1}{c(\varsigma)} \E\Bracket{\psi^{\prime\prime}(X_1\T\theta) \cdot X_1X_1\T}.
    \end{align*}
    The interchange of differentiation and integration above is justified by Corollary \ref{cor:regularity}, which is applicable since $\ell$ is $\nu_\ell$-smooth and has a uniformly bounded gradient. 

    To establish a Loewner order lower bound on $\nabla_\theta^2 \mathcal L_P(\theta)$, define where 
    \begin{align*}
        \mathcal E 
        := \curly{\abs{\theta\T X_1} \leq M_{\upsilon,\zeta,\Theta}}
        \quad \text{where} \quad 
        M_{\upsilon,\zeta,\Theta} 
        := \sqrt{2\upsilon^2 R_\Theta^2 \log\Parent{128 \cdot \frac{\upsilon^4}{\zeta^2}}}. 
    \end{align*}
    Herein, the logarithm is positive as $\zeta \leq \var{v\T X} \leq \upsilon^2$ for all unit length $v \in \mathcal S^{d-1}$ and thus $\frac{\upsilon^2}{\zeta} \geq 1$. By the symmetry and positivity of $\psi^{\prime\prime}$ and as it is non-increasing, for all $v \in \R^d$,
    \begin{align*}
        v\T \nabla_\theta^2 \mathcal L_P(\theta) v
        &= \frac{1}{c(\varsigma)} \E\Bracket{\psi^{\prime\prime}(X_1\T\theta) (v\T X_1)^2}
        = \frac{1}{c(\varsigma)} \E\Bracket{\psi^{\prime\prime}(\abs{X_1\T\theta}) (v\T X_1)^2} \\
        &= \frac{1}{c(\varsigma)} \Parent{\E\Bracket{\psi^{\prime\prime}(\abs{X_1\T\theta}) (v\T X_1)^2 \indicator{\mathcal E}} + \E\Bracket{\psi^{\prime\prime}(\abs{X_1\T\theta}) (v\T X_1)^2 \indicator{\mathcal E^c}}} \\
        &\geq \frac{1}{c(\varsigma)}  \psi^{\prime\prime}(M_{\upsilon,\zeta,\Theta}) \cdot \E\Bracket{(v\T X_1)^2 \indicator{\mathcal E}}.
    \end{align*}
    Using that $\E\Bracket{X_1X_1\T} \succeq \zeta I_d$, by definition of the Loewner order,
    \begin{align*}
        \E\Bracket{(v\T X_1)^2 \indicator{\mathcal E}}
        &= \E\Bracket{(v\T X_1)^2} - \E\Bracket{(v\T X_1)^2 \indicator{\mathcal E^c}} \\
        &= v\T \E\Bracket{X_1 X_1\T} v - \E\Bracket{(v\T X_1)^2 \indicator{\mathcal E^c}} \\
        &\geq \zeta \twonorm{v}^2 - \E\Bracket{(v\T X_1)^2 \indicator{\mathcal E^c}} \\
        &\geq \frac{\zeta}{2} \twonorm{v}^2.
    \end{align*}
    It remains to justify the last inequality above. Because $X_1$ is a zero-mean $\upsilon^2$-sub-Gaussian vector, $v\T X_1$ and $\theta\T X_1$ are both zero-mean sub-Gaussian random variables in $\R$. Combined with Cauchy Schwarz, this yields the following control on the 4-th moment and tail:
    \begin{align*}
        \E\Bracket{(v\T X_1)^2 \indicator{\mathcal E^c}}
        &\leq \sqrt{\E\Bracket{(v\T X_1)^4}} \sqrt{\P\Parent{\mathcal E^c}} \\
        &= \sqrt{\E\Bracket{(v\T X_1)^4}} \sqrt{\P\Parent{\abs{\theta\T X_1} > M}} \\
        &\leq \sqrt{16 \upsilon^4 \twonorm{v}^4} \sqrt{2\exp\Parent{-\frac{M_{\upsilon,\zeta,\Theta}^2}{2\upsilon^2 \twonorm{\theta}^2}}} \\
        &\overset{(i)} \leq 4 \upsilon^2 \twonorm{v}^2 \sqrt{2\exp\Parent{-\frac{M_{\upsilon,\zeta,\Theta}^2}{2\upsilon^2 R_\Theta^2}}} \\
        &\overset{(ii)} \leq \frac{\zeta}{2} \twonorm{v}^2.
    \end{align*}
    Here, $(i)$ holds as $\twonorm{\theta} \leq R_\Theta$ and $(ii)$ by definition of $M_{\upsilon,\zeta,\Theta}$. Collecting terms, for all $v \in \R^d$
    \begin{align*}
        v\T \nabla_\theta^2 \mathcal L_P(\theta) v
        \geq \frac{\zeta}{2c(\varsigma)} \cdot \psi^{\prime\prime}(M_{\upsilon,\zeta,\Theta}) \cdot \twonorm{v}^2
        \nabla_\theta^2 \mathcal L_P(\theta) \succeq \frac{\zeta}{2c(\varsigma)} \cdot \psi^{\prime\prime}(M_{\upsilon,\zeta,\Theta}) \cdot I_d. 
    \end{align*}
    Therefore, $\mathcal L_P$ is  $\mu$-strongly convex with $\mu := \frac{\zeta}{2c(\varsigma)} \psi^{\prime\prime}(M_{\upsilon,\zeta,\Theta}) $. The constants so far are
    \begin{align*}
        L := \parent{C_\psi^\prime + R_{\mathcal Y}} \frac{R_{\mathcal X}}{c(\varsigma)} \sqrt{d}, \quad 
        \nu_\ell := \frac{1}{4c(\varsigma)}C_\psi^{\prime\prime} R_{\mathcal X}^2d
        \quad \text{and} \quad 
        \mu := \frac{\zeta}{2c(\varsigma)} \psi^{\prime\prime}(M_{\upsilon,\zeta,\Theta})
    \end{align*}
    Hence, Theorem \ref{thm:esterrorinexp} with $b = 1$, $\beta = 4$ and $\tau = \lceil \frac{\nu_\ell}{\mu} +2\rceil \leq 4\frac{\nu_\ell}{\mu} \asymp \frac{\nu_\ell}{\mu}$ yields 
    \begin{align*}
        \E\Bracket{\twonorm{\theta_n - \theta^\star}^2}
        &\lesssim \E\Bracket{\twonorm{\theta_0 - \theta^\star}^2} \cdot \Parent{\frac{\tau}{n}}^4 + \frac{L^2}{\mu^2 n} + \frac{\alpha dL^2}{\mu^2\epsilon n^2} \\
        &\lesssim R_\Theta^2 \cdot \frac{\nu_\ell}{\mu n} + \frac{L^2}{\mu^2} \Parent{\frac{1}{n} + \frac{\alpha d}{\epsilon n^2}} \\
        &\lesssim R_\Theta^2 \cdot \frac{C_\psi^{\prime\prime} R_{\mathcal X}^2}{\zeta \psi^{\prime\prime}(M_{\upsilon,\zeta,\Theta})} \Parent{\frac{d}{n}} + \frac{\parent{C_\psi^\prime + R_{\mathcal Y}}^2 R_{\mathcal X}^2}{\zeta^2 \psi^{\prime\prime}(M_{\upsilon,\zeta,\Theta})^2} \Parent{\frac{d}{n} + \frac{\alpha d^2}{\epsilon n^2}}.
    \end{align*}
    The last inequality holds by plugging in the values for $L, \nu_\ell$ and $\mu$. Note here, that the assumption that $\tau < n$ in Theorem \ref{thm:esterrorinexp} is fulfilled as we impose the constraint
    \begin{align*}
        n
        \gtrsim \frac{C_\psi^{\prime\prime} R_{\mathcal X}^2 d}{\zeta \psi^{\prime\prime}(M_{\upsilon,\zeta,\Theta})}.
    \end{align*}
    As $\mathcal L_P$ is strongly convex, $\theta^\star$ is its unique minimizer. By Lemma \ref{lem:recoveryPRM} the true parameter $\vartheta$ also minimizes $\mathcal L_P$. We can therefore conclude that $\theta^\star = \vartheta$, which recovers the statement. 
\end{proof}

\begin{corollary}
    \label{cor:sigmaprimelb}
    Let $\sigma$ be the sigmoid function. Then, $\sigma^\prime(x) \geq \frac 1 4 e^{-\abs{x}}$ for all $x \in \R$. 
\end{corollary}

\begin{proof}
    By definition of the sigmoid function for all $x \in \R$:  
    \begin{align*}
        \sigma^\prime(x) 
        = \sigma(x) (1-\sigma(x))
        = \sigma(x) \sigma(-x)
        = \frac{1}{1+e^{-x}} \frac{1}{1+e^{x}}
        = \frac{1}{2 + e^x + e^{-x}}
        \geq \frac 1 4 e^{-\abs{x}}.
    \end{align*}
    The inequality above holds by a case distinction. Using that $e^x \geq 1$ and $e^{-x} \leq 1$ for $x \geq 0$ and $e^{-x} \geq 1$ and $e^x \leq 1$ for $x < 0$ we may justify the following: 
    \begin{align*}
        2 + e^x + e^{-x}
        \leq 
        \begin{cases} 
            3 + e^x
            \leq 4 e^x
            = 4 e^{\abs{x}} & \text{if $x \geq 0$, }\\
            3 + e^{-x}
            \leq 4 e^{-x}
            = 4 e^{\abs{x}} & \text{if $x < 0$}. 
        \end{cases}
    \end{align*}
\end{proof}

\section{Auxiliary Results}

\begin{lemma}
    \label{lem:unrollrec}
    For $\alpha_k \in [0,1]$ and $A_k, B_k \in \R$, assume the following holds for all $k \in [K]$: 
    \begin{align*}
        \gamma_k \cdot A_{k-1} + B_k
        \leq A_k
        \leq \alpha_k \cdot A_{k-1} + B_k. 
    \end{align*}
    Then, unrolling the recursions implied by the lower and upper bound up to $K$ yields 
    \begin{align*}
        \Parent{\prod_{j=1}^K \gamma_j} \cdot A_0 + \sum_{j=1}^K B_j \prod_{i=j+1}^K \gamma_i
        \leq A_K
        \leq \Parent{\prod_{j=1}^K \alpha_j} \cdot A_0 + \sum_{j=1}^K B_j \prod_{i=j+1}^K \alpha_i. 
    \end{align*}
\end{lemma}

\begin{proof}
    For $k = K$ the assumed upper bound yields 
    \begin{align*}
        A_K
        &\leq \alpha_K \cdot A_{K-1} + B_K \\
        &\leq \alpha_K \cdot \Parent{\alpha_{K-1} \cdot A_{K-2} + B_{K-1}} + B_K \\
        &= \alpha_K \cdot \alpha_{K-1} \cdot A_{K-2} + \alpha_K \cdot B_{K-1} + B_K \\
        &\leq \alpha_K \cdot \alpha_{K-1} \cdot \Parent{\alpha_{K-2} \cdot A_{K-3} + B_{K-2}} + \alpha_K \cdot B_{K-1} + B_K \\
        &= \alpha_K \cdot \alpha_{K-1} \cdot \alpha_{K-2} \cdot A_{K-3} + \alpha_K \cdot \alpha_{K-1} \cdot B_{K-2} + \alpha_K \cdot B_{K-1} + B_K. 
    \end{align*}
    Unrolling this down to $k = 0$ and using product and sum notation yields the statement   
    \begin{align*}
        A_K
        \leq \Parent{\prod_{j=1}^K \alpha_j} \cdot A_0 + \sum_{j=1}^K B_j \prod_{i=j+1}^K \alpha_i. 
    \end{align*}
    The proof for the $\geq$-case is analogous and instead uses the lower bounds on $A_k$. 
\end{proof}

\begin{corollary}
    \label{cor:product}
    Let $j, \beta, \tau \in \N$ where $\tau \geq \beta-1$. Then, the product simplifies into $\beta$ terms: 
    \begin{align*}
        \prod_{i=j+1}^n \Parent{1-\frac{\beta}{i+\tau}}
        = \frac{(j+\tau - (\beta-1)) \dotsm (j+\tau)}{(n+\tau - (\beta-1)) \dotsm (n+\tau)}.
    \end{align*}
\end{corollary}

\begin{proof}
    By explicitly writing out the statement's \lhs{}, we obtain
    \begin{align*}
        \prod_{i=j+1}^n \Parent{1-\frac{\beta}{i+\tau}}
        &= \prod_{i=j+1}^n \Parent{\frac{i+\tau-\beta}{i+\tau}} \\
        &= \Parent{\frac{j+1+\tau-\beta}{j+1+\tau}}\Parent{\frac{j+2+\tau-\beta}{j+2+\tau}} \dotsm \Parent{\frac{n+\tau-\beta}{n+\tau}} \\
        &= \Parent{\frac{j+\tau-(\beta-1)}{j+\tau+1}}\Parent{\frac{j+\tau-(\beta-2)}{j+\tau+2}} \dotsm \Parent{\frac{n+\tau-\beta}{n+\tau}} \\
        &= \frac{(j+\tau-(\beta-1))(j+\tau-(\beta-2))\dotsm(j+\tau)}{(n+\tau-(\beta-1))(n+\tau-(\beta-2)) \dotsm (n+\tau)}. 
    \end{align*}
\end{proof}

\begin{corollary}
    \label{cor:sumSGD}
    Let $\eta_k=\frac{\mu+\nu}{2\mu\nu}\frac{\beta}{k+\tau}$ for $k\in[K]$ and $\tau \geq \beta-1 \geq 2$. Then, 
    \begin{align*}
        \sum_{k=1}^K \eta_k^2 \Pi_{k,K}^{\mu,\nu}
        \leq \frac{\beta^2}{K} \Parent{\frac{\mu+\nu}{2\mu\nu}}^2.
    \end{align*}
\end{corollary}

\begin{proof}
    By Corollary \ref{cor:product} we have
    \begin{align*}
        \sum_{k=1}^K \eta_k^2 \Pi_{k,K}^{\mu,\nu}
        &= \sum_{k=1}^K \eta_k^2 \frac{(k+\tau - (\beta-1))(k+\tau - (\beta-2)) \dotsm (k+\tau)}{(K+\tau - (\beta-1))(K+\tau - (\beta-2)) \dotsm (K+\tau)} \\
        &= \Parent{\frac{\mu+\nu}{2\mu\nu}}^2 \beta^2 \sum_{k=1}^K \frac{1}{(k+\tau)^2}\frac{(k+\tau - (\beta-1))(k+\tau - (\beta-2)) \dotsm (k+\tau)}{(K+\tau - (\beta-1))(K+\tau - (\beta-2)) \dotsm (K+\tau)} \\
        &\leq \Parent{\frac{\mu+\nu}{2\mu\nu}}^2 \beta^2 \sum_{k=1}^K \frac{(k+\tau - (\beta-1))(k+\tau - (\beta-2)) \dotsm (k+\tau-2)}{(K+\tau - (\beta-1))(K+\tau - (\beta-2)) \dotsm (K+\tau)} \\
        &\leq \Parent{\frac{\mu+\nu}{2\mu\nu}}^2 \beta^2 \sum_{k=1}^K \frac{1}{(K+\tau-1)(K+\tau)} \\
        &\leq \Parent{\frac{\mu+\nu}{2\mu\nu}}^2 \frac{\beta^2}{K}. 
    \end{align*}
\end{proof}

\begin{corollary}
    \label{cor:sumDP}
    Take $\eta_k=\frac{\mu+\nu}{2\mu\nu}\frac{\beta}{k+\tau}$ for all $k\in[K]$ and choose $\curly{\sigma^2_k} \in [0,\infty)^K$ as in Remark \ref{rem:conversion}. Assume $\cardinality{B_k} = b$ for all $k \in [K-1]$ and $\cardinality{B_K} \geq b$. Then, for $\tau \geq \beta-1 \geq 3$, 
    \begin{align*}
        \sum_{k=1}^K  \sigma_k^2 \eta_k^2 \Pi_{k,K}^{\mu,\nu}
        \leq \frac{3\alpha L^2\beta^2}{\epsilon b^2K^2} \Parent{\frac{\mu+\nu}{2\mu\nu}}^2.
    \end{align*}
\end{corollary}

\begin{proof}
    We start this proof with two tedious computations that we will use later. 

    First, for any $k \in [K]$ the following computations hold
    \begin{align*}
        \sum_{k=1}^K \frac{\beta^2}{(k+\tau)^3} \Pi_{k,K}^{\mu,\nu}
        &=\sum_{k=1}^K \frac{\beta^2}{(k+\tau)^3} \frac{(k+\tau - (\beta-1)) \dotsm (k+\tau)}{(K+\tau - (\beta-1)) \dotsm (K+\tau)} \\
        &\leq \beta^2\sum_{k=1}^K \frac{(k+\tau - (\beta-1)) \dotsm (k+\tau-3)}{(K+\tau - (\beta-1)) \dotsm (K+\tau)} \\
        &\leq \beta^2\sum_{k=1}^K \frac{1}{(K+\tau-2)(K+\tau-1)(K+\tau)} \\
        &\leq \beta^2 \frac{1}{(K+\tau-1)(K+\tau)} \\
        &\leq \frac{\beta^2}{K^2}.
    \end{align*}
    Second, for any $k \in [K-1]$ we have
    \begin{align*}
        \eta_k^2-\eta_{k+1}^2
        &= \Parent{\frac{\mu+\nu}{2\mu\nu}}^2 \beta^2 \Parent{\frac{1}{(k+\tau)^2}-\frac{1}{(k+\tau+1)^2}} \\
        &= \Parent{\frac{\mu+\nu}{2\mu\nu}}^2 \beta^2 \Parent{\frac{(k+\tau+1)^2-(k+\tau)^2}{(k+\tau)^2(k+\tau+1)^2}} \\
        &= \Parent{\frac{\mu+\nu}{2\mu\nu}}^2 \beta^2 \Parent{\frac{k^2+\tau^2+1+2k\tau+2k+2\tau-k^2-\tau^2-2k\tau}{(k+\tau)^2(k+\tau+1)^2}} \\
        &= \Parent{\frac{\mu+\nu}{2\mu\nu}}^2 \beta^2 \Parent{\frac{1+2k+2\tau}{(k+\tau)^2(k+\tau+1)^2}} \\
        &\leq \Parent{\frac{\mu+\nu}{2\mu\nu}}^2 \frac{2\beta^2}{(k+\tau)^2(k+\tau+1)} \\
        &\leq \Parent{\frac{\mu+\nu}{2\mu\nu}}^2 \frac{2\beta^2}{(k+\tau)^3}.
    \end{align*}
    By definition of the $\Pi_{k,K}^{\mu,\nu}$, it holds that
    \begin{align*}
        \Pi_{k,K}^{\mu,\nu} - \Pi_{k-1,K}^{\mu,\nu}
        = \Pi_{k,K}^{\mu,\nu} - \Parent{1-\frac{2\eta_k\mu\nu}{\mu+\nu}}\Pi_{k,K}^{\mu,\nu}
        = \frac{2\eta_k\mu\nu}{\mu+\nu} \Pi_{k,K}^{\mu,\nu}
        \geq 0
    \end{align*}
    Combining this with $\cardinality{B_k} \geq b$ for all $k \in [K]$ which holds by assumption and $\Pi_{K,K}^{\mu,\nu} = 1$,
    \begin{align*}
        \sum_{k=1}^{K-1} &\Parent{\frac{\eta_k^2}{\cardinality{B_k}^2} - \frac{\eta_{k+1}^2}{\cardinality{B_{k+1}}^2}} \Pi_{k,K}^{\mu,\nu} + \frac{\eta_K^2}{\cardinality{B_K}^2} \Pi_{K,K}^{\mu,\nu} \\
        &= \frac{\eta_1^2}{\cardinality{B_1}} \Pi_{1,K}^{\mu,\nu} + \sum_{k=2}^K \frac{\eta_k^2}{\cardinality{B_k}} \parent{\Pi_{k,K}^{\mu,\nu}- \Pi_{k-1,K}^{\mu,\nu}} \\
        &\leq \frac{1}{b} \Parent{\eta_1^2 \Pi_{1,K}^{\mu,\nu} + \sum_{k=2}^K \eta_k^2 \parent{\Pi_{k,K}^{\mu,\nu}- \Pi_{k-1,K}^{\mu,\nu}}} \\
        &= \frac{1}{b} \Parent{\sum_{k=1}^{K-1} \Parent{\eta_k^2 - \eta_{k+1}^2} \Pi_{k,K}^{\mu,\nu} + \eta_K^2}
    \end{align*}
    Given a noise schedule as in Remark \ref{rem:conversion}, 
    \begin{align*}
        \sum_{k=1}^K \sigma_k^2 \eta_k^2 \Pi_{k,K}^{\mu,\nu}
        &= \frac{\alpha L^2}{\epsilon} \Parent{\sum_{k=1}^{K-1} \Parent{\frac{\eta_k^2}{\cardinality{B_k}^2} - \frac{\eta_{k+1}^2}{\cardinality{B_{k+1}}^2}} \Pi_{k,K}^{\mu,\nu} + \frac{\eta_K^2}{\cardinality{B_K}^2} \Pi_{K,K}^{\mu,\nu}} \\
        &\overset{(i)} \leq \frac{\alpha L^2}{\epsilon b^2} \Parent{\sum_{k=1}^{K-1} \Parent{\eta_k^2 - \eta_{k+1}^2} \Pi_{k,K}^{\mu,\nu} + \eta_K^2} \\
        &\leq \frac{\alpha L^2}{\epsilon b^2} \Parent{\frac{\mu+\nu}{2\mu\nu}}^2\Parent{\sum_{k=1}^{K-1} \frac{2\beta^2}{(k+\tau)^3} \Pi_{k,K}^{\mu,\nu} + \frac{\beta^2}{(K+\tau)^2}} \\
        &\leq \frac{\alpha L^2}{\epsilon b^2} \Parent{\frac{\mu+\nu}{2\mu\nu}}^2 \Parent{\frac{2\beta^2}{K^2} + \frac{\beta^2}{(K+\tau)^2}} \\
        &\overset{(ii)} \leq \frac{3\alpha L^2\beta^2}{\epsilon b^2 K^2} \Parent{\frac{\mu+\nu}{2\mu\nu}}^2. 
    \end{align*}
    Here, $(i)$ was shown above and $(ii)$ holds because $\frac{1}{(K+\tau)^2} \leq \frac{1}{K^2}$. 
\end{proof}

\begin{corollary}
    \label{cor:prodBias}
    Take $\eta_k=\frac{\mu+\nu}{2\mu\nu}\frac{\beta}{k+\tau}$ for all $k\in[K]$. For $\tau \geq \beta-1 \geq 0$, we have 
    \begin{align*}
        \Pi_{0,K}^{\mu,\nu}
        \leq \Parent{\frac{\tau}{K}}^{\beta}. 
    \end{align*}
\end{corollary}

\begin{proof}
    By Corollary \ref{cor:product} we have
    \begin{align*}
        \Pi_{0,K}^{\mu,\nu}
        &= \frac{(\tau - (\beta-1))(\tau - (\beta-2)) \dotsm \tau}{(K+\tau - (\beta-1))(K+\tau - (\beta-2)) \dotsm (K+\tau)}
        \leq \Parent{\frac{\tau}{K}}^{\beta}. 
    \end{align*}
    Above, the inequality holds because $\tau \geq \beta-l \geq 0$ for all $l \in [\beta]$.  
\end{proof}

\begin{corollary}
    \label{cor:prodsimplification}
    Let $0 < \mu \leq \nu < \infty$ and $\{\eta_k\} \in [0,\infty)^K$. Then, 
    \begin{align*}
        \Parent{1-\frac{2\eta_k \mu\nu}{\mu+\nu}} \Pi_{l,k-1}^{\mu,\nu} \Pi_{k,K}^{\mu,\nu}
        = \Pi_{l,K}^{\mu,\nu}. 
    \end{align*}
\end{corollary}

\begin{proof}
    By definition of $\Pi_{l,n}^{\mu,\nu}$, for all $l \in [k-1]$, 
    \begin{align*}
        \Parent{1-\frac{2\eta_k \mu\nu}{\mu+\nu}} \Pi_{l,k-1}^{\mu,\nu} \Pi_{k,K}^{\mu,\nu}
        &= \Parent{1-\frac{2\eta_k \mu\nu}{\mu+\nu}}\prod_{j=l+1}^{k-1} \Parent{1-\frac{2\eta_j\mu\nu}{\mu+\nu}} \prod_{j=k+1}^K \Parent{1-\frac{2\eta_j\mu\nu}{\mu+\nu}} \\
        &= \prod_{j=l+1}^K \Parent{1-\frac{2\eta_j\mu\nu}{\mu+\nu}}
        = \Pi_{l,K}^{\mu,\nu}. 
    \end{align*}
\end{proof}

\end{document}